\documentclass[12pt]{article}
\usepackage{color}

\usepackage{amsmath}
\usepackage{amssymb}
\usepackage{amsbsy}
\usepackage{graphicx}
\usepackage{cite}

\usepackage{amsfonts}
\newtheorem{lem}{Lemma}[section]%
\newtheorem{theorem}[lem]{Theorem}%
\newtheorem{prop}[lem]{Proposition}%

\def\a{\alpha} \def\b{\beta} \def\g{\gamma} \def\d{\delta} 
    
 \def\ld{\lambda} 

\def\G{\Gamma}

 \def\lg{\langle} \def\rg{\rangle}
\def\nd{\mathrel{\bigm|\kern-.7em/}}

\def\f{\noindent}

\def\Aut{\hbox{\rm Aut\,}}

\def\Cay{\hbox{\rm Cay}}
\def\Sym{\hbox{\rm Sym}}
\def\BiCay{\hbox{\rm BiCay}}

\def\mod{\hbox{\rm mod }}

\def\Core{\hbox{\rm Core}}

\def\C{\hbox{\rm C}}
\def\cl{\hbox{\rm CL}}
\def\mcl{\hbox{\rm MCL}}

\def\R{{\mathcal R}}
\def\mb{{\mathcal B}}
\def\demo{{\bf Proof}\hskip10pt}

\def\mz{{\mathbb Z}}
\def\qed{\hskip10pt $\Box$\vspace{3mm}}

\begin{document}

\title{Trivalent dihedrants and bi-dihedrants}
\author {{\sc Mi-Mi Zhang} , {\sc Jin-Xin Zhou}
 \\{\small{Department of Mathematics, Beijing Jiaotong University, Beijing, 100044, China}}}
\date{}
\maketitle

\begin{abstract}

{\rm\small \noindent{~~~~
A Cayley (resp. bi-Cayley) graph on a dihedral group is called a {\em dihedrant} (resp. {\em bi-dihedrant}).
In 2000, a classification of trivalent arc-transitive dihedrants was given by Maru\v si\v c and Pisanski, and several years later, trivalent non-arc-transitive dihedrants of order $4p$ or $8p$ $(p$ a prime) were classified by Feng et al. As a generalization of these results, our first result presents a classification of trivalent non-arc-transitive dihedrants. Using this, a complete classification of trivalent vertex-transitive non-Cayley bi-dihedrants is given, thus completing the study of trivalent bi-dihedrants initiated in our previous paper [Discrete Math. 340 (2017) 1757--1772]. As a by-product, we generalize a theorem in [The Electronic Journal of Combinatorics 19 (2012) $\#$P53].}
\bigskip

}
\end{abstract}

\section{Introduction}

In this paper we describe an investigation of trivalent Cayley graphs on dihedral groups as well as vertex-transitive trivalent bi-Cayley graphs over dihedral groups. To be brief, we shall say that a Cayley (resp. bi-Cayley) graph on a dihedral group a {\em dihedrant} (resp. {\em bi-dihedrant}).

Cayley graphs are usually defined in the following way. Given a finite group $G$ and an inverse closed subset $S\subseteq G\setminus\{1\}$, the {\em Cayley graph} $\Cay(G,S)$ on $G$ with respect to $S$ is a graph with vertex set $G$ and edge set $\{\{g,sg\}\mid g\in G,s\in S\}$. For any $g\in G$, $R(g)$ is the permutation of $G$ defined by $R(g): x\mapsto xg$ for $x\in G$.  Set $R(G):=\{R(g)\ |\ g\in G\}$. It is well-known that $R(G)$ is a subgroup of $\Aut(\Cay(G,S))$. We say that the {Cayley graph} $\Cay(G,S)$ is {\em normal} if $R(G)$ is normal in $\Aut(\Cay(G,S))$ (see \cite{X1}).

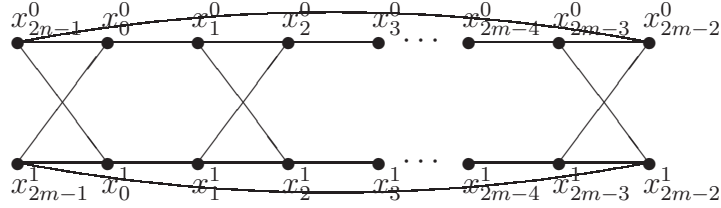
\begin{figure}[ht]
\begin{center}
\unitlength 4mm
\begin{picture}(25,8)

\put(2, 6){\line(1, 0){12}}\put(14.8, 5.8){$\cdots$}\put(17, 6){\line(1, 0){6}}
\put(2, 2){\line(3, 4){3}}\put(2, 6){\line(3, -4){3}} \put(8, 2){\line(3, 4){3}}\put(8, 6){\line(3, -4){3}}  \put(20, 2){\line(3, 4){3}}\put(20, 6){\line(3, -4){3}}
\put(2, 2){\line(1, 0){12}}\put(14.8, 1.8){$\cdots$}\put(17, 2){\line(1, 0){6}}

\put(1.75, 5.7){$\bullet$}\put(4.75, 5.7){$\bullet$}\put(7.75, 5.7){$\bullet$}\put(10.75, 5.7){$\bullet$}\put(13.75, 5.7){$\bullet$}
\put(16.75, 5.7){$\bullet$}\put(19.75, 5.7){$\bullet$}\put(22.75, 5.7){$\bullet$}

\put(1.75, 1.7){$\bullet$}\put(4.75, 1.7){$\bullet$}\put(7.75, 1.7){$\bullet$}\put(10.75, 1.7){$\bullet$}\put(13.75, 1.7){$\bullet$}
\put(16.75, 1.7){$\bullet$}\put(19.75, 1.7){$\bullet$}\put(22.75, 1.7){$\bullet$}

\put(1.8, 6.5){$x_{2n-1}^{0}$}\put(4.8, 6.5){$x_{0}^{0}$}\put(7.8, 6.5){$x_{1}^{0}$}\put(10.8, 6.5){$x_{2}^{0}$}\put(13.8, 6.5){$x_{3}^{0}$}
 \put(16.8, 6.5){$x_{2m-4}^{0}$}\put(19.8, 6.5){$x_{2m-3}^{0}$}\put(22.8, 6.5){$x_{2m-2}^{0}$}

\put(1.8, 1){$x_{2m-1}^{1}$}\put(4.8, 1){$x_{0}^{1}$}\put(7.8, 1){$x_{1}^{1}$}\put(10.8, 1){$x_{2}^{1}$}\put(13.8, 1){$x_{3}^{1}$}
\put(16.8, 1){$x_{2m-4}^{1}$}\put(19.8, 1){$x_{2m-3}^{1}$}\put(22.8, 1){$x_{2m-2}^{1}$}

\bezier{600}(2, 6)(12.5,8)(23,6)
\bezier{600}(2, 2)(12.5,0)(23,2)

\end{picture}
\end{center}\vspace{-.5cm}
\caption{The cross ladder graph $\cl_{4m}$} \label{Fig-1}
\end{figure}

In 2000, Maru\v si\v c and Pisanski~\cite{MP} initiated the study of automorphisms of dihedrants, and they gave a classification of trivalent arc-transitive dihedrants. Following this work, highly symmetrical dihedrants have been extensively studied, and one of the remarkable achievements is the complete classification of 2-arc-transitive dihedrants (see \cite{DMM,M2008}). In contrast, however, relatively little is known about the automorphisms of non-arc-transitive dihedrants. In \cite{A-G-S-normal-dihedrant}, the authors proved that every trivalent non-arc-transitive dihedrant is normal. However, this is not true. {There exist non-arc-transitive and non-normal dihedrants.} Actually, in \cite{ZhouCX-Feng,ZM}, the automorphism groups of trivalent dihedrants of order $4p$ and $8p$ are determined for each prime $p$, and the result reveals that every non-arc-transitive trivalent dihedrant of order $4p$ or $8p$ is either a normal Cayley graph, or isomorphic to the so-called cross ladder graph. For an integer $m\geq 2$, the {\em cross ladder graph}, denoted by $\cl_{4m}$, is a trivalent graph of order $4m$ with vertex set $V_0\cup V_1\cup \ldots V_{2m-2}\cup V_{2m-1}$, where $V_i=\{x_i^0,x_i^1\}$,
and edge set $\{\{x_{2i}^r, x_{2i+1}^r\}, \{x^r_{2i+1}, x^s_{2i+2}\}\mid i\in\mz_m, r,s\in\mz_2\}$ (see Fig.~$\ref{Fig-1}$ for $\cl_{4m}$).
It is worth mentioning that the cross ladder graph plays an important role in the study of automorphisms of trivalent graphs (see, for example, \cite{Cameron-Spiga2006,ZM,ZFZ}). Motivated by the above mentioned facts, we shall focus on trivalent non-arc-transitive dihedrants. Our first theorem generalizes the results in \cite{ZhouCX-Feng,ZM} to all trivalent dihedrants.

\begin{theorem}\label{not-at-dihedrant}
Let $\Sigma=\Cay(H, S)$ be a connected trivalent Cayley graph, where $H=\lg a, b\ |\ a^n=b^2=1, bab=a^{-1}\rg (n\geq 3)$. If $\Sigma$ is non-arc-transitive and non-normal, then $n$ is even and $\G\cong \cl_{4\cdot\frac{n}{2}}$ and $S^\a=\{b, ba, ba^{\frac{n}{2}}\}$ for some $\a\in\Aut(H)$.
\end{theorem}

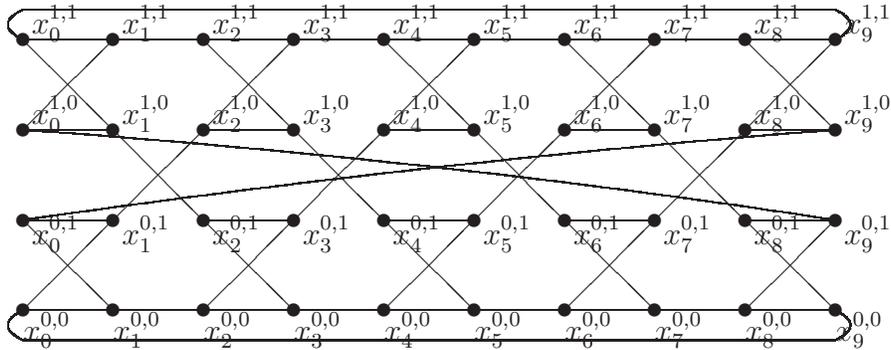
\begin{figure}[ht]
\begin{center}
\unitlength 4mm
\begin{picture}(29,11)
\put(1, 1){\circle*{0.4}} \put(1,1){\line(1,0){27}}\put(1,1){\line(1,1){3}}\put(1,0){\line(1,0){27}}

\qbezier(1, 1)(0, 0.5)(1, 0)\qbezier(28, 1)(29, 0.5)(28, 0)

\put(4, 1){\circle*{0.4}}

\put(7, 1){\circle*{0.4}}\put(7,1){\line(1,1){3}}

\put(10, 1){\circle*{0.4}}

\put(13, 1){\circle*{0.4}}\put(13,1){\line(1,1){3}}

\put(16, 1){\circle*{0.4}}

\put(19, 1){\circle*{0.4}}\put(19,1){\line(1,1){3}}

\put(22, 1){\circle*{0.4}}

\put(25, 1){\circle*{0.4}}\put(25,1){\line(1,1){3}}

\put(28, 1){\circle*{0.4}}



\put(1, 4){\circle*{0.4}}\put(1,4){\line(1,0){3}}\put(1,4){\line(1,-1){3}}

\put(4, 4){\circle*{0.4}}\put(4,4){\line(1,1){3}}

\put(7, 4){\circle*{0.4}}\put(7,4){\line(1,-1){3}}\put(7,4){\line(1,0){3}}

\put(10, 4){\circle*{0.4}}\put(10,4){\line(1,1){3}}

\put(13, 4){\circle*{0.4}}\put(13,4){\line(1,-1){3}}\put(13,4){\line(1,0){3}}

\put(16, 4){\circle*{0.4}}\put(16,4){\line(1,1){3}}

\put(19, 4){\circle*{0.4}}\put(19,4){\line(1,-1){3}}\put(19,4){\line(1,0){3}}

\put(22, 4){\circle*{0.4}}\put(22,4){\line(1,1){3}}

\put(25, 4){\circle*{0.4}}\put(25,4){\line(1,-1){3}}\put(25,4){\line(1,0){3}}

\put(28, 4){\circle*{0.4}}\qbezier(28, 4)(15, 6)(1, 7)


\put(1, 7){\circle*{0.4}}\put(1,7){\line(1,0){3}}\put(1,7){\line(1,1){3}}

\put(4, 7){\circle*{0.4}}\put(4,7){\line(1,-1){3}}

\put(7, 7){\circle*{0.4}}\put(7,7){\line(1,0){3}}\put(7,7){\line(1,1){3}}

\put(10, 7){\circle*{0.4}}\put(10,7){\line(1,-1){3}}

\put(13, 7){\circle*{0.4}}\put(13,7){\line(1,0){3}}\put(13,7){\line(1,1){3}}

\put(16, 7){\circle*{0.4}}\put(16,7){\line(1,-1){3}}

\put(19, 7){\circle*{0.4}}\put(19,7){\line(1,0){3}}\put(19,7){\line(1,1){3}}

\put(22, 7){\circle*{0.4}}\put(22,7){\line(1,-1){3}}

\put(25, 7){\circle*{0.4}}\put(25,7){\line(1,0){3}}\put(25,7){\line(1,1){3}}

\put(28, 7){\circle*{0.4}}\qbezier(28, 7)(15, 6)(1, 4)



\put(1, 10){\circle*{0.4}}\put(1,10){\line(1,0){27}}\put(1,10){\line(1,-1){3}}\put(1,11){\line(1,0){27}}

\qbezier(1, 10)(0, 10.5)(1, 11)\qbezier(28, 10)(29, 10.5)(28, 11)

\put(4, 10){\circle*{0.4}}

\put(7, 10){\circle*{0.4}}\put(7,10){\line(1,-1){3}}

\put(10, 10){\circle*{0.4}}

\put(13, 10){\circle*{0.4}}\put(13,10){\line(1,-1){3}}

\put(16, 10){\circle*{0.4}}

\put(19, 10){\circle*{0.4}}\put(19,10){\line(1,-1){3}}

\put(22, 10){\circle*{0.4}}

\put(25, 10){\circle*{0.4}}\put(25,10){\line(1,-1){3}}

\put(28, 10){\circle*{0.4}}



\put(1, 0){$x_0^{0,0}$}
\put(4, 0){$x_1^{0,0}$}
\put(7, 0){$x_2^{0,0}$}
\put(10, 0){$x_3^{0,0}$}
\put(13, 0){$x_4^{0,0}$}
\put(16, 0){$x_5^{0,0}$}
\put(19, 0){$x_6^{0,0}$}
\put(22, 0){$x_7^{0,0}$}
\put(25, 0){$x_8^{0,0}$}
\put(28, 0){$x_9^{0,0}$}

\put(1.3, 10.2){$x_0^{1,1}$}
\put(4.3, 10.2){$x_1^{1,1}$}
\put(7.3, 10.2){$x_2^{1,1}$}
\put(10.3, 10.2){$x_3^{1,1}$}
\put(13.3, 10.2){$x_4^{1,1}$}
\put(16.3, 10.2){$x_5^{1,1}$}
\put(19.3, 10.2){$x_6^{1,1}$}
\put(22.3, 10.2){$x_7^{1,1}$}
\put(25.3, 10.2){$x_8^{1,1}$}
\put(28.3, 10.2){$x_9^{1,1}$}

\put(1.3, 7.2){$x_0^{1,0}$}
\put(4.3, 7.2){$x_1^{1,0}$}
\put(7.3, 7.2){$x_2^{1,0}$}
\put(10.3, 7.2){$x_3^{1,0}$}
\put(13.3, 7.2){$x_4^{1,0}$}
\put(16.3, 7.2){$x_5^{1,0}$}
\put(19.3, 7.2){$x_6^{1,0}$}
\put(22.3, 7.2){$x_7^{1,0}$}
\put(25.3, 7.2){$x_8^{1,0}$}
\put(28.3, 7.2){$x_9^{1,0}$}

\put(1.3, 3.2){$x_0^{0,1}$}
\put(4.3, 3.2){$x_1^{0,1}$}
\put(7.3, 3.2){$x_2^{0,1}$}
\put(10.3, 3.2){$x_3^{0,1}$}
\put(13.3, 3.2){$x_4^{0,1}$}
\put(16.3, 3.2){$x_5^{0,1}$}
\put(19.3, 3.2){$x_6^{0,1}$}
\put(22.3, 3.2){$x_7^{0,1}$}
\put(25.3, 3.2){$x_8^{0,1}$}
\put(28.3, 3.2){$x_9^{0,1}$}
\end{picture}
\end{center}\vspace{-.5cm}
\caption{The multi-cross ladder graph $\mcl_{20, 2}$} \label{Fig-2}
\end{figure}

Recall that for an integer $m\geq2$, the {cross ladder graph} $\cl_{4m}$ has vertex set $V_0\cup V_1\cup \ldots V_{2m-2}\cup V_{2m-1}$, where $V_i=\{x_i^0,x_i^1\}$. The {\em multi-cross ladder graph}, denoted by $\mcl_{4m, 2}$, is the graph obtained from $\cl_{4m}$ by blowing up each
vertex $x_i^r$ of $\cl_{4m}$ into two vertices $x_i^{r,0}$ and $x_i^{r,1}$. The edge set is
$\{\{x^{r,s}_{2i}, x^{r,t}_{2i+1}\}, \{x^{r,s}_{2i+1}, x^{s,r}_{2i+2}\}\mid i\in\mz_m, r,s,t\in\mz_2\}$ (see Fig.~$\ref{Fig-2}$ for $\mcl_{20, 2}$).

Note that the multi-cross ladder graph $\mcl_{4m, 2}$ is just the graph given in \cite[Definition~7]{ZhouF}. From \cite[Proposition~3.3]{DMMN} we know that every $\mcl_{4m, 2}$ is vertex-transitive. However, not all
multi-cross ladder graphs are Cayley graphs. Actually, in \cite[Theorem~9]{ZhouF}, it is proved that $\mcl_{4p, 2}$ is a vertex-transitive non-Cayley graph for each prime $p>7$.  Our second theorem generalizes this result to all multi-cross ladder graphs.

\begin{theorem}\label{EX(m,2)-Cayley}
The multi-cross ladder graph $\mcl_{4m, 2}$ is a Cayley graph if and only if either $m$ is even, or $m$ is odd and $3\mid m$.
\end{theorem}

Both of the above two theorems are crucial in attacking the problem of classification of trivalent vertex-transitive non-Cayley bi-dihedrants. Before proceeding, we give some background to this topic, and set some notation.

Let $R,L$ and $S$ be subsets of a group $H$ such that $R=R^{-1}$, $L=L^{-1}$ and $R\cup L$ does not contain the identity element of $H$. The {\em bi-Cayley graph} BiCay$(H, R, L, S)$ over $H$ relative to $R, L, S$ is a graph having vertex set the union of the {\em right part} $H_0=\{h_0\ |\ h\in H\}$ and the {\em left part} $H_1=\{h_1\ |\ h\in H\}$, and edge set the union of the {\em right edges} $\{\{h_0,g_0\}\ |\ gh^{-1}\in R\}$, the {\em left edges} $\{\{h_1,g_1\}\ |\ gh^{-1}\in L\}$ and the {\em spokes} $\{\{h_0,g_1\}\ |\ gh^{-1}\in S\}$. If $|R|=|L|=s$, then $\BiCay(H,~R,~L,~S)$ is said to be an {\em $s$-type bi-Cayley graph}.

In \cite{ZZ} we initiated a program to investigate the automorphism groups of the trivalent vertex-transitive bi-dihedrants. This was partially motivated by the following facts. As one of the most important finite graphs, the Petersen graph is a bi-circulant, but it is not a Cayley graph.
Note that a {\em bi-circulant} is a bi-Cayley graph over a cyclic group.
The Petersen graph is the initial member of a family of graphs $P(n,t)$, known now as
the {\em generalized Petersen graphs} (see~\cite{Watkins}), which can be also constructed as bi-circulants.
Let $n\geq 3$, $1\leq t<n/2$ and set $H=\lg a\rg\cong\mz_n$.
The {generalized Petersen graph $P(n,t)$} is isomorphic to
the bi-circulant ${\rm BiCay}(H,~\{a, a^{-1}\},~\{a^t, a^{-t}\},~\{1\})$.
The complete classification of vertex-transitive generalized Petersen graphs has been worked out
in~\cite{FGW,NS}. Latter, this was generalized by Maru\v si\v c et al. in \cite{MP,Pisanski} where all trivalent vertex-transitive
bi-circulants were classified in~\cite{MP,Pisanski}, and more recently, all trivalent vertex-transitive bi-Cayley graphs over abelian groups were classified in~\cite{ZF}. The characterization of trivalent vertex-transitive bi-dihedrants is the next natural step.

Another motivation for us to consider trivalent vertex-transitive bi-dihedrants comes from the work in \cite{Potocnik1280}. By checking the census of trivalent vertex-transitive graphs of order up to 1000 in \cite{Potocnik1280}, we find out that there are $981$ non-Cayley graphs, and among these graphs, $233$ graphs are non-Cayley bi-dihedrants. This may suggest bi-dihedrants form an important class of trivalent vertex-transitive non-Cayley graphs.

In \cite{ZZ}, we gave a classification of trivalent arc-transitive bi-dihedrants, and we also proved that every trivalent vertex-transitive $0$- or $1$-type bi-dihedrant is a Cayley graph, and gave a classification of trivalent vertex-transitive non-Cayley bi-dihedrants of order $4n$ with $n$ odd. The goal of this paper is to complete the classification of trivalent vertex-transitive non-Cayley bi-dihedrants.

%



Before stating the main result, we need the following concepts.
For a bi-Cayley graph $\Gamma=\BiCay(H,~R,~L,~S)$ over a group $H$, we can assume that the identity $1$ of $H$
is in $S$ (see Proposition~\ref{properties}~(2)). The triple $(R, L, S)$ of three subsets $R, L, S$ of a group $H$ is called {\em bi-Cayley triple} if $R=R^{-1},  L=L^{-1}, $ and $1\in S$. Two bi-Cayley triples $(R, L, S)$ and $(R', L', S')$ of a group $H$ are said to be
{\em equivalent}, denoted by $(R, L, S)\equiv(R', L', S')$, if either $(R', L', S')=(R, L, S)^\a$ or $(R', L', S')=(L, R, S^{-1})^\a$
for some automorphism $\a$ of $H$. The bi-Cayley graphs corresponding to
two equivalent bi-Cayley triples of the same group are isomorphic (see Proposition~\ref{properties}~(3)-(4)).

\begin{theorem}\label{all-non-cay}
Let $\G=\BiCay(R,L,S)$ be a trivalent vertex-transitive bi-dihedrant where $H=\lg a,b \ |\ a^{n}=b^{2}=1,bab=a^{-1}\rg$ is a dihedral group. Then either $\G$ is a Cayley graph or one of the following occurs:
\begin{enumerate}
\item [$(1)$]  $(R, L, S)\equiv(\{b,~ba^{\ell+1}\}, \{ba,~ba^{\ell^{2}+\ell+1}\}, \{1\})$, where $n\geq 5$,
$\ell^{3}+\ell^{2}+\ell+1\equiv0~(\mod n)$, $\ell^{2}\not\equiv 1~(\mod n)$.

\item [$(2)$]  $(R, L, S)\equiv(\{ba^{-\ell},~ba^{\ell}\}, \{a,~a^{-1}\}, \{1\})$, where $n=2m$ and $\ell^{2}\equiv-1~(\mod m)$. Furthermore, $\G$ is also a bi-Cayley graph over an abelian group $\mathbb{Z}_{n}\times \mathbb{Z}_{2}$.

\item [$(3)$] $(R,L,S)\equiv (\{b,ba\},\{b,ba^{2m}\},\{1\})$, where $n=2(2m+1)$, $m\not\equiv 1~(\mod 3)$, and the corresponding graph is isomorphic the multi-cross ladder graph $\mcl_{4m,2}$.

\item [$(4)$] $(R,L,S)\equiv (\{b,ba\},\{ba^{24\ell},ba^{12\ell-1}\},\{1\})$, where $n=48\ell$ and $\ell\geq 1$.
\end{enumerate}
Moreover, all of the graphs arising from $(1)$-$(4)$ are vertex-transitive non-Cayley.

\end{theorem}

\section{Preliminaries}
All groups considered in this paper are finite, and all graphs are finite, connected, simple and undirected.
For the group-theoretic and graph-theoretic terminology not defined here we refer the
reader to \cite{BMBook,WI}.

\subsection{Definitions and notations}
For a positive integer, let $\mathbb{Z}_{n}$ be the cyclic group of order $n$
and $\mathbb{Z}_{n}^{*}$ be the multiplicative group of $\mathbb{Z}_{n}$ consisting of numbers coprime to $n$.
For two groups $M$ and $N$, $N\rtimes M$ denotes a semidirect product of $N$ by $M$.
For a subgroup $H$ of a group $G$, denote $C_{G}(H)$  the centralizer of $H$ in $G$ and by $N_{G}(H)$ the normalizer of $H$ of $G$.
Let $G$ be a permutation group on a set $\Omega$ and $\a \in \Omega$.
Denote by $G_{\a}$ the stabilizer of $\a$ in $G$.
We say that $G$ is {\em semiregular} on $\Omega$ if $G_{\a}=1$ for every $\a\in\Omega$
and {\em regular} if $G$ is transitive and semiregular.

For a finite, simple and undirected graph $\G$, we use $V(\G)$, $E(\G)$, $A(\G)$, $\Aut(\G)$ to denote its vertex set, edge set, arc set and full automorphism group, respectively. 
For any subset $B$ of $V(\G)$, the subgraph of $\G$ induced by $B$
will be denoted by $\G[B]$. For any $v\in V(\G)$ and a positive integer $i$ no more than the diameter of $\G$, denote by $\G_i(v)$ be the set of vertices at distance $i$ from $v$. Clearly, $\G_1(v)$ is just the neighborhood of $v$. We shall often abuse the notation by using $\G(v)$ to replace $\G_1(v)$.

A graph $\G$ is said to be {\em vertex-transitive}, and {\em arc-transitive} (or {\em symmetric}) if $\Aut(\G)$ acts transitively on $V(\G)$ and $A(\G)$, respectively. Let $\G$ be a connected vertex-transitive
graph, and let $G\leq \Aut(\G)$ be vertex-transitive on $\G$. For a
$G$-invariant partition $\mathcal{B}$ of $V(\G)$, the {\em quotient graph}
$\G_{\mathcal{B}}$ is defined as the graph with vertex set $\mathcal{B}$ such that, for
any two different vertices $B, C\in \mathcal{B}$, $B$ is adjacent to $C$ if and only if
there exist $u\in B$ and $v\in C$ which are adjacent in $\G$. Let $N$
be a normal subgroup of $G$. Then the set $\mathcal{B}$ of orbits of $N$ in
$V(\G)$ is a $G$-invariant partition of $V(\G)$. In this case, the
symbol $\G_{\mathcal{B}}$ will be replaced by $\G_N$. The original graph $\G$ is
said to be a {\em $N$-cover} of $\G_N$ if $\G$ and $\G_N$ have the same valency.

\subsection{Cayley graphs}

Let $\G=\Cay(G,S)$ be a Cayley graph on $G$ with respect to $S$. Then $\G$ is vertex-transitive due to $R(G)\leq\Aut(\G)$. In
general, we have the following proposition.

\begin{prop}{\rm \cite[Lemma~16.3]{B}}\label{IsCay}
A vertex-transitive graph $\G$ is isomorphic to a Cayley
graph on a group $G$ if and only if its automorphism group has a
subgroup isomorphic to $G$, acting regularly on the vertex set of
$\G$.
\end{prop}

In 1981, Godsil~\cite{Godsil1981} proved that the normalizer of $R(G)$ in $\Aut(\Cay(G,S))$ is $R(G)\rtimes\Aut(G,S)$, where $\Aut(G,S)$ is the group of automorphisms of $G$ fixing the set $S$ set-wise. This result has been successfully used in characterizing various families of Cayley graphs $\Cay(G,S)$ such that $R(G)=\Aut(\Cay(G,S))$ (see, for example, \cite{Godsil1981,Godsil1983}).
Recall that a Cayley graph $\Cay(G, S)$ is
said to be {\em normal} if $R(G)$ is normal in $\Aut(\Cay(G, S))$ (see~\cite{X1}).

\begin{prop} {\rm\cite[Proposition~1.5]{X1}}  \label{stabi}
The Cayley graph $\G=\Cay(G,S)$ is normal if and only if
$A_1=\Aut(G,S)$, where $A_1$ is the stabilizer of the identity $1$
of $G$ in $\Aut(\G)$.
\end{prop}

\subsection{Basic properties of bi-Cayley graphs}

In this subsection, we let $\G$ be a connected bi-Cayley graph $\BiCay(H,R,L,S)$
over a group $H$. It is easy to prove some basic properties of such a $\G$, as in \cite[Lemma~3.1]{ZF}.

\begin{prop}\label{properties}
 The following hold.
\begin{enumerate}
\item[$(1)$] $H$ is generated by $R\cup L\cup S$.
\item[$(2)$] Up to graph isomorphism, $S$ can be chosen to contain the identity of $H$.
\item[$(3)$] For any automorphism $\alpha$ of $H$, $\BiCay(H,~R,~L,~S)\cong \BiCay(H,~R^{\alpha},~L^{\alpha},~S^{\alpha})$.
\item[$(4)$] $\BiCay(H,~R,~L,~S) \cong \BiCay(H,~L,~R,~S^{-1})$.
\end{enumerate}
\end{prop}

Next, we collect several results about the automorphisms of bi-Cayley graph $\G={\rm BiCay}(H,~R,~L,~S)$.
For each $g\in H$, define a permutation as follows:
\begin{equation}\label{1}
\R(g): h_{i}\mapsto (hg)_{i},~~~\forall i\in \mathbb{Z}_{2},~h\in H.
\end{equation}
Set $\R(H)=\{\R(g)\ |\ g\in H\}$. Then $\R(H)$ is a semiregular subgroup of $\Aut(\G)$ with $H_0$ and $H_1$ as its two orbits.

For an automorphism $\alpha$ of $H$ and $x, y, g\in H$, define two permutations of $V(\Gamma)=H_{0}\cup H_{1}$ as follows:
\begin{equation} \label{2}
\begin{aligned}
&\delta_{\a,x,y}:~ h_{0}\mapsto(xh^{\a})_{1},~ h_{1}\mapsto(yh^{\a})_{0},~\forall h\in H,\\
&~~\sigma_{\a,g}:~ h_{0}\mapsto(h^{\alpha})_{0},~ h_{1}\mapsto(gh^{\a})_{1},~\forall h\in H.\\
\end{aligned}
\end{equation}

Set
\begin{equation} \label{3}
\begin{aligned}
&I=\{ \delta_{\a,x,y}~|~\alpha \in \Aut(H)~s.t.~R^{\alpha}=x^{-1}Lx,~L^{\alpha}=y^{-1}Ry,~S^{\alpha}=y^{-1}S^{-1}x\},\\
&F=\{\sigma_{\alpha,g}~|~\alpha \in \Aut(H)~s.t.~R^{\alpha}=R,~L^{\alpha}=g^{-1}Lg,~S^{\alpha}=g^{-1}S\}.\\
\end{aligned}
\end{equation}

\begin{prop}{\rm\cite[Theorem~1.1]{ZF-auto}}\label{normalizer}
Let $\Gamma=\BiCay(H,R,L,S)$ be a connected bi-Cayley graph over the group $H$.
Then $N_{\Aut(\Gamma)}(\R(H))=\R(H)\rtimes F$ if $I=\emptyset$ and
$N_{\Aut(\Gamma)}(\R(H))=\R(H)\langle F,~\delta_{\a,x,y}\rangle$ if
$I\neq\emptyset$ and $\delta_{\a,x,y}\in I$. Furthermore, for any $\delta_{\a,x,y}\in I$, we have the following:
\begin{enumerate}
\item[$(1)$] $\langle \R(H),~\delta_{\a,x,y}\rangle$ acts transitively on $V(\Gamma)$;
\item[$(2)$] if $\a$ has order $2$ and $x=y=1$, then $\Gamma$ is isomorphic to the Cayley graph Cay($\bar{H},~R\cup \alpha S$), where $\bar{H}=H\rtimes \langle\a\rangle$.
\end{enumerate}
\end{prop}

\section{Cross ladder graphs}

The goal of this section is to prove Theorem~\ref{not-at-dihedrant}.\medskip

\f{\bf Proof of Theorem~\ref{not-at-dihedrant}}\ Suppose that $\Sigma=\Cay(H, S)$ is a connected trivalent Cayley graph which is neither normal nor arc-transitive, where $H=\lg a, b\ |\ a^n=b^2=1, bab=a^{-1}\rg (n\geq 3)$. Then $S$ is a generating subset of $H$ and $|S|=3$. So $S$ must contain an involution of $H$ outside $\lg a\rg$. As $\Aut(H)$ is transitive on the coset $b\lg a\rg$, we may assume that $S=\{b, x, y\}$ for $x, y\in H\setminus\lg b\rg$.

Suppose first that $x$ is not an involution. Then we must have $y=x^{-1}$. Since $S$ generates $H$, one has $\lg a\rg=\lg x\rg$, and so $bxb=x^{-1}$. Then there exists an automorphism of $H$ sending $b, x$ to $b, a$ respectively. So we may assume that $S=\{b, a, a^{-1}\}$.
Now it is easy to check that $\Sigma$ is isomorphic to the generalized Petersen graph { $P(n, 1)$}. Since $\Sigma$ is not arc-transitive, by \cite{FGW,NS}, we have $|\Aut(\Sigma)|=2|H|$, and so $\Sigma$ would be a normal Cayley graph of $H$, a contradiction.

Therefore, both $x$ and $y$ must be involutions. Suppose that $x\in\lg a\rg$. Then $n$ is even and $x=a^{n/2}$.
Again since $S$ generates $H$, one has $y=ba^{j}$, where $1\leq j\leq n-1$ and either $(j, n)=1$ or $(j, n)=2$ and $\frac{n}{2}$ is odd.
Note that the subgroup of $\Aut(H)$ fixing $b$ is transitive on the set of generators of $\lg a\rg$ and that
$\lg a^{n/2}\rg$ is the center of $H$. There exists $\a\in\Aut(H)$ such that \[S^\a=\{b, ba, a^{\frac{n}{2}}\}\ {\rm or}\ \{b, ba^{2}, a^{\frac{n}{2}}\}.\]
Without loss of generality, we may assume that $S=\{b, ba, a^{\frac{n}{2}}\}$ or $\{b, ba^{2}, a^{\frac{n}{2}}\}.$
If $S=\{b, ba^2, a^{\frac{n}{2}}\}$, {we shall prove that $\Sigma\cong P(n,1)$. Note that the generalized Petersen graph $P(n,1)$ has vertex set $\{u_{i},v_{i}\ |\ i\in \mz_{n}\}$ and edge set $\{\{u_{i},u_{i+1}\}, \{v_{i},v_{i+1}\}, \{u_{i},v_{i}\}\ |\ i\in \mz_{n}\}$.
Define a map from $V(\Sigma)$ to $V(P(n,1))$ as follows:
$$\begin{array}{lll}
\varphi: & a^{2i}\mapsto u_{2i}, &a^{2i+\frac{n}{2}}\mapsto v_{2i},\\
      & ba^{2i}\mapsto u_{2i-1},&ba^{2i+\frac{n}{2}}\mapsto v_{2i-1},
\end{array}$$ where $0\leq i\leq \frac{n}{2}-1$.
It is easy to see that $\varphi$ is an isomorphism form $\Sigma$ to $P(n, 1)$}. Since $\Sigma$ is not arc-transitive, by \cite{FGW,NS}, we have $|\Aut(\Sigma)|=2|H|$, and so $\Sigma$ would be a normal Cayley graph of $H$, a contradiction.
If $S=\{b, ba, a^{\frac{n}{2}}\}$, then $\Sigma$ has a connected subgraph $\Sigma_1=\Cay(H, \{b, ba\})$ which is a cycle of length $2n$, and $\Sigma$ is just the graph obtained from $\Sigma_1$ by adding a $1$-factor such that each vertex $g$ of $\Sigma_1$ is adjacent to its antipodal vertex $a^{\frac{n}{2}}g$. Then $R(H)\rtimes\mz_2\cong \Aut(\Sigma_1)\leq\Aut(\Sigma)$, and then since $\Sigma$ is assumed to be not arc-transitive, $\Aut(\Sigma)$ will fix the $1$-factor $\{\{g, a^{\frac{n}{2}}g\}\ |\ g\in H\}$ setwise. This implies that $\Aut(\Sigma)\leq\Aut(\Sigma_1)$ and so $\Aut(\Sigma)=\Aut(\Sigma_1)$. Consequently, we have $\Sigma$ is a  normal Cayley graph of $H$, a contradiction.

Similarly, we have $y\notin\lg a\rg$. Then we may assume that $x=ba^{i}$ and $y=ba^{j}$ for some $1\leq i, j\leq n-1$ and $i\neq j$.
Then $S=\{b, ba^{i}, ba^{j}\}\subseteq b\lg a\rg$. This implies that $\Sigma$ is a bipartite graph with $\lg a\rg$ and $b\lg a\rg$ as its two partition sets. Since $\Sigma$ is not arc-transitive, $\Aut(\Sigma)_1$ is intransitive on the neighbourhood $S$ of $1$, and since $\Sigma$ is not a normal Cayley graph of $H$, there exists a unique element, say $s\in S$, such that $\Aut(\Sigma)_1=\Aut(\Sigma)_s$.
Considering the fact that $\Aut(H)$ is transitive on $b\lg a\rg$, without loss of generality, we may assume that $\Aut(\Sigma)_1=\Aut(\Sigma)_b$ and $\Aut(\Sigma)_1$ swaps $ba^i$ and $ba^j$. Then for any $h\in H$, we have \[\Aut(\Sigma)_h=(\Aut(\Sigma)_1)^{R(h)}=(\Aut(\Sigma)_b)^{R(h)}=\Aut(\Sigma)_{bh}.\]

Direct computation shows that
$$
\begin{array}{l}
  \Sigma_{2}(1)=\{a^{-i}, a^{-j}, a^{i}, a^{i-j}, a^{j}, a^{j-i} \}, \\
  \Sigma_{3}(1)=\{ba^{-i}, ba^{j-i}, ba^{-j}, ba^{i-j}, ba^{2i}, ba^{j+i}, ba^{2i-j}, ba^{2j}, ba^{2j-i}\}.
   \end{array}
$$

Let $\Aut(\Sigma)_1^{*}$ be the kernel of $\Aut(\Sigma)_1$ acting on $S$. Take an $\a\in \Aut(\Sigma)_1^{*}$. Then $\a$ fixes every element in $S$.
As $\Aut(\Sigma)_h=\Aut(\Sigma)_{bh}$ for any $h\in H$, $\a$ will fix $b(ba^{i})=a^{i}$ and $b(ba^{j})=a^{j}$.
Note that $\Sigma(ba^{i})\setminus \{1, a^{i}\}=\{a^{i-j}\}$ and $\Sigma(ba^{j})\setminus\{1, a^{j}\}=\{a^{j-i}\}$.
Then $\a$ also fixes $a^{i-j}$ and $a^{j-i}$, and then $\a$ also fixes $ba^{i-j}$ and $ba^{j-i}$.

If $|\Sigma_{2}(1)|=6$, then it is easy to check that $a^{-i}$ is the unique common neighbor of $b$ and $ba^{j-i}$. So $\a$ also fixes $a^{-i}$. Now one can see that $\a$ fixes every vertex in $\Sigma_2(1)$. If $|\Sigma_2(1)|<6$ and either $|\Sigma_1(b)\cap\Sigma_1(ba^i)|>1$ or $|\Sigma_1(b)\cap\Sigma_1(ba^j)|>1$, then $\a$ also fixes every vertex in $\Sigma_2(1)$.
In the above two cases, by the connectedness and vertex-transitivity of $\Sigma$, $\a$ would fix all vertices of $\Sigma$, implying that $\a=1$. Hence, $\Aut(\Sigma)_1^{*}=1$ and $\Aut(\Sigma)_1\cong \mz_2$. This forces that $\Sigma$ is a normal Cayley graph of $H$, a contradiction.


Thus, we have $|\Sigma_2(1)|<6$ and $|\Sigma_1(b)\cap\Sigma_1(ba^i)|=|\Sigma_1(b)\cap\Sigma_1(ba^j)|=1$. This implies that $\Sigma_1(ba^i)\cap\Sigma_1(ba^j)=\{1, a^{i-j}\}=\{1, a^{j-i}\}$, and so $a^{i-j}=a^{j-i}$. It follows that $a^{i-j}$ is an involution, and hence $n$ is even and $a^{i-j}=a^{n/2}$. So $S=\{b, ba^i, ba^{i+n/2}\}$. As $S$ generates $H$, one has $\lg a^i, a^{n/2}\rg=\lg a\rg$.
So either $(i, n)=1$ or $(i, n)=2$ and $\frac{n}{2}$ is odd.
Note that the subgroup of $\Aut(H)$ fixing $b$ is transitive on the set of generators of $\lg a\rg$ and that
$\lg a^{n/2}\rg$ is the center of $H$. There exists $\a\in\Aut(H)$ such that \[S^\a=\{b, ba, ba^{1+\frac{n}{2}}\}\ {\rm or}\ \{b, ba^{2}, ba^{2+\frac{n}{2}}\}.\]
Let $\b$ be the automorphism of $H$ induced by the map $a\mapsto a^{-1}, b\mapsto ba$.
Then
\[\{b, ba, ba^{1+\frac{n}{2}}\}^\b=\{b, ba, ba^{\frac{n}{2}}\}, {\rm and}\ \{b, ba^{2}, ba^{2+\frac{n}{2}}\}^\b=\{b, ba^{2}, ba^{\frac{n}{2}}\}. \]
If $\frac{n}{2}$ is odd, then the map $\eta: a\mapsto a^{2+\frac{n}{2}}, b\mapsto ba^{\frac{n}{2}}$ induces an automorphism of $H$, and $\{b, ba, ba^{\frac{n}{2}}\}^{\eta}=\{b, ba^{2}, ba^{\frac{n}{2}}\}$.
So there always exists $\g\in\Aut(H)$ such that $S^\g=\{b, ba, ba^{\frac{n}{2}}\}$, completing the proof of the first part of our theorem.

Finally, we shall prove $\Sigma\cong \cl_{4\cdot\frac{n}{2}}$. Without loss of generality, assume that $S=\{b, ba, ba^{\frac{n}{2}}\}$.
Recall that $V(\cl_{4\cdot\frac{n}{2}})=\{x_{i}^{r}\ |\ i\in\mz_{2n},r\in\mz_{2}\}$ and $E(\cl_{4\cdot\frac{n}{2}})=\{\{x_{i}^{r}, x_{i+1}^{r}\},$ $\{x_{2i}^{r}, x_{2i+1}^{r+1}\},\ |\ i\in\mz_{2n},r\in\mz_{2}\}$.
Let $\phi$ be a map from $V(\Sigma)$ to $V(\cl_{4\cdot\frac{n}{2}})$ as following:
$$\begin{array}{lll}
\phi: & a^{i}\mapsto x_{2i}^{0}, &a^{i+\frac{n}{2}}\mapsto x_{2i}^{1},\\
      & ba^{j}\mapsto x_{2j-1}^{0},&ba^{j+\frac{n}{2}}\mapsto x_{2j-1}^{1},
\end{array}$$ where $0\leq i\leq \frac{n}{2}-1$ and $1\leq j\leq \frac{n}{2}$.
It is easy to check that $\phi$ is an isomorphism from $\Sigma$ and $X(\cl_{4\cdot\frac{n}{2}})$, as desired.
\hfill\qed

\section{Multi-cross ladder graphs}

The goal of this section is to prove Theorem~\ref{EX(m,2)-Cayley}. We first show that each $\mcl_{4m, 2}$ is a bi-Cayley graph.

\begin{lem}\label{EX(m,2)}
The multi-cross ladder graph $\mcl_{4m, 2}$ is isomorphic to the bi-Cayley graph $\BiCay(H,\{c,ca\},\{ca,ca^{2}b\},\{1\})$, where
\[H=\lg a,b,c\ |\ a^{m}=b^{2}=c^{2}=1,a^b=a,a^{c}=a^{-1},b^{c}=b\rg.\]
\end{lem}

\f\demo For convenience, let $\G$ be the bi-Cayley graph given in our lemma, and let $X=\mcl_{4m, 2}$. Let $\phi$ be a map from $V(X)$ to $V(\G)$ defined by the following rule:
$$\begin{array}{lllll}
\phi: &x_{2t}^{1,1}\mapsto (a^{t})_{0}, &x_{2t+1}^{1,1}\mapsto (ca^{t+1})_{0},&x_{2t}^{1,0}\mapsto (ca^{t+1})_{1}, &x_{2t+1}^{1,0}\mapsto (a^{t})_{1},\\
~\\
      &x_{2t}^{0,1}\mapsto (ca^{t+1}b)_{1}, &x_{2t+1}^{0,1}\mapsto (a^{t}b)_{1}, &x_{2t}^{0,0}\mapsto (a^{t}b)_{0}, &x_{2t+1}^{0,0}\mapsto (ca^{t+1}b)_{0},
\end{array}$$ where $t\in\mz_{m}$.

It is easy to see that $\phi$ is a bijection. Furthermore, for any $t\in\mz_{m}$, we have
$$
\begin{array}{l}
X(x_{2t}^{1,1})^\phi=\{x_{2t+1}^{1,0}, x_{2t+1}^{1,1},
x_{2t-1}^{1,1}\}^{\phi}=\{(a^{t})_{1},(ca^{t+1})_{0},(ca^{t})_{0}\}=\G((a^{t})_{0}),\\
X(x_{2t+1}^{1,1})^\phi=\{x_{2t}^{1,0}, x_{2t+2}^{1,1}, x_{2t}^{1,1}\}^{\phi}=\{(ca^{t+1})_{1},(a^{t+1})_{0},(a^{t})_{0}\}=\G((ca^{t+1})_{0}),\\
X(x_{2t}^{1,0})^\phi=\{x_{2t+1}^{1,1}, x_{2t+1}^{1,0}, x_{2t-1}^{0,1}\}^{\phi}=\{(ca^{t+1})_{0},(a^{t})_{1},(a^{t-1}b)_{1}\}=\G((ca^{t+1})_{1}),\\
X(x_{2t+1}^{1,0})^\phi=\{x_{2t}^{1,1}, x_{2t}^{1,0}, x_{2t+2}^{0,1}\}^{\phi}=\{(a^{t})_{0},(ca^{t+1})_{1},(ca^{t+2}b)_{1}\}=\G((a^{t})_{1}),\\
X(x_{2t}^{0,1})^\phi=\{x_{2t+1}^{0,0}, x_{2t+1}^{0,1}, x_{2t-1}^{1,0}\}^{\phi}=\{(ca^{t+1}b)_{0},(a^{t}b)_{1},(a^{t-1})_{1}\}=\G((ca^{t+1}b)_{1}),\\
X(x_{2t+1}^{0,1})^\phi=\{x_{2t}^{0,0}, x_{2t}^{0,1}, x_{2t+2}^{1,0}\}^{\phi}=\{(a^{t}b)_{0},(ca^{t+1}b)_{1},(ca^{t+2})_{1}\}=\G((a^{t}b)_{1}),\\
X(x_{2t}^{0,0})^\phi=\{x_{2t+1}^{0,1}, x_{2t+1}^{0,0}, x_{2t-1}^{0,0}\}^{\phi}=\{(a^{t}b)_{1},(ca^{t+1}b)_{0},(ca^{t}b)_{0}\}=\G((a^{t}b)_{0}),\\
X(x_{2t+1}^{0,0})^\phi=\{x_{2t}^{0,1}, x_{2t}^{0,0}, x_{2t+2}^{0,0}\}^{\phi}=\{(ca^{t+1}b)_{1},(a^{t}b)_{0},(a^{t+1}b)_{0}\}=\G((ca^{t+1}b)_{0}).
\end{array}
$$
This shows that $\phi$ preserves the adjacency of the graphs, and so it is an isomorphism from $X$ to $\G$.
\hfill\qed

\f{\bf Remark~1}\  Let $m$ be odd, let $e=ab$ and $f=ca$. Then the group given in Lemma~\ref{EX(m,2)} has the following presentation:
 \[H=\lg e, f \ |\ e^{2m}=f^{2}=1, e^{f}=e^{-1}\rg.\]
Clearly, in this case, $H$ is a dihedral group.
Furthermore, the corresponding bi-Cayley graph given in Lemma~\ref{EX(m,2)} will be
\[\BiCay(H, \{f, fe\}, \{f, fe^{m-1}\}, \{1\}).\]

\f{\bf Proof of Theorem~\ref{EX(m,2)-Cayley}}\ By Lemma~\ref{EX(m,2)}, we may let $\G=\mcl_{4m, 2}$ be just the bi-Cayley graph
$\BiCay(H, R, L, S),$ where
\[\begin{array}{l}
H=\lg a,b,c\ |\ a^{m}=b^{2}=c^{2}=1,a^b=a,a^{c}=a^{-1},b^{c}=b\rg,\\
R=\{c,ca\}, L=\{ca,ca^{2}b\}, S=\{1\}.
\end{array}\]

We first prove the sufficiency. Assume first that $m$ is even. Then the map
\[a\mapsto ab, b\mapsto b, c\mapsto cb\]
induces an automorphism, say $\a$ of $H$ of order $2$.
Furthermore, $R^\a=\{c, ca\}^\a=caLca$, $L^\a=\{ca,ca^{2}b\}^\a=caRca$ and $S^\a=\{1\}^\a=ca\{1\}ca=S^{-1}$.
By Proposition~\ref{normalizer}, $\d_{\a, ca, ca}\in\Aut(\G)$ and $\R(H)\rtimes\lg\d_{\a, ca, ca}\rg$ acts regularly on $V(\G)$.
Consequently, by Proposition~\ref{IsCay}, $\G$ is a Cayley graph.

Assume now that $m$ is odd and $3\mid m$. In this case, we shall use the bi-Cayley presentation for $\G$ as in Remark~5.1, that is, \[\G=\BiCay(H,\{f,fe\},\{f,fe^{m-1}\},\{1\}),\] where
\[H=\lg e, f\ |\ e^{2m}=f^{2}=1,e^{f}=e^{-1}\rg.\]
Let $\b$ be a permutation of $V(\G)$ defined as following:
$$\begin{array}{lll}
\b: &(f^{i}e^{3t+1})_{i}\leftrightarrow (f^{i}e^{m+3t+1})_{i},&(f^{i+1}e^{3t+1})_{i}\leftrightarrow (f^{i}e^{m+3t+1})_{i+1}, \\
    &(f^{i+1}e^{3t+2})_{i}\leftrightarrow (f^{i+1}e^{m+3t+2})_{i},&(f^{i}e^{3t+2})_{i}\leftrightarrow (f^{i+1}e^{m+3t+2})_{i+1},\\
    &(e^{3t})_{i}\leftrightarrow (fe^{3t})_{i+1}, & (e^{m+3t})_{i}\leftrightarrow (fe^{m+3t})_{i+1},\\
\end{array}$$ where $t\in\mz_{\frac{m}{3}}$ and $i\in\mz_{2}$.
It is easy to check that $\b$ is an automorphism of $\G$ of order $2$. Furthermore, $\R(e),\R(f)$ and $\b$
satisfy the following relations:
$$\begin{array}{l}\label{relation}
\R(e)^{2m}=\R(f)^{2}=\b^{2}=1,~\R(f)^{-1}\R(e)\R(f)=\R(e)^{-1},~\R(f)^{-1}\b\R(f)=\b,\\
\R(e)^{6}\b=\b\R(e)^{6},~\R(e)^{2}\b=\b \R(e)^{4}\b \R(e)^{-2}.
\end{array}$$
Let $G=\lg \R(e^{2}),\R(f), \b\rg$ and $P=\lg \R(e^{2}), \b\rg$.
Then $\R(f)\notin P$ and $G=P\lg \R(f)\rg$. Since $\R(e)^{6}\b=\b\R(e)^{6}$, we have $\R(e^{6})\in Z(P)$. Since $\R(e)^{2}\b=\b \R(e)^{4}\b \R(e)^{-2}$, it follows that
$$(\R(e)^{2}\b)^{3}=\R(e)^{2}\b[\b \R(e)^{4}\b \R(e)^{-2}]\R(e)^{2}\b=\R(e^{6}).$$
Let $N=\lg \R(e^{6})\rg$.
Clearly, $N$ is a normal subgroup of $G$. Furthermore,
\[P/N=\lg \R(e^{2})N, \b N\ |\ \R(e^{2})^{3}N=\b^{2}N=(\R(e^{2})\b)^{3}N=N\rg\cong A_{4}.\] Therefore, $|P|=4m$ and $|G|\leq8m$.

Let $$\begin{array}{ll}
\Delta_{00}=\{x_{0}\ |\ x\in \lg e^{2},f\rg\}, &\Delta_{10}=\{(ex)_{0}\ |\ x\in \lg e^{2},f\rg\},\\
\Delta_{01}=\{x_{1}\ |\ x\in \lg e^{2},f\rg\}, &\Delta_{11}=\{(ex)_{1}\ |\ x\in \lg e^{2},f\rg\}.
\end{array}$$
Then $\Delta_{ij}$'s $(i,j\in\mz_{2})$ are four orbits of $\lg \R(e^{2}),\R(f)\rg$.
Moreover,
$$1_{0}^{\b\R(f)}=1_{1}\in \Delta_{01},~e_{0}^{\b}=(e^{m+1})_{0}\in \Delta_{00},~e_{1}^{\b}=(fe^{m+1})_{0}\in \Delta_{00}.$$
This implies that $G$ is transitive on $V(\G)$.
Hence, $|G|=8m$ and so $G$ is regular on $V(\G)$, and by Proposition~\ref{IsCay}, $\G$ is a Cayley graph.

To prove the necessity, it suffices to prove that if $m$ is odd and $3\nmid m$, then $\G$ is a non-Cayley graph.
In this case, we shall use the original definition of $\G=\mcl_{4m,2}$. Suppose that $m$ is odd and $3\nmid m$.
We already know from \cite[Proposition~3.3]{DMMN} that $\G$ is vertex-transitive. Let $A=\Aut(\G)$.
For $m=5 $ or $7$, using Magma \cite{MAGMA}, $\G$ is a non-Cayley graph. In what follows, we assume that $m\geq 11$.

For each $j\in\mz_{m}$, $\C_{j}^{0}=(x_{2j}^{0,0},x_{2j+1}^{0,0},x_{2j}^{0,1},x_{2j+1}^{0,1})$ and $\C_{j}^{1}=(x_{2j}^{1,1},x_{2j+1}^{1,1},x_{2j}^{1,0},x_{2j+1}^{1,0})$ are two $4$-cycles.
Set $\mathcal{F}=\{\C_{j}^{i}\ |\ i\in\mz_{2},j\in\mz_{m}\}$. From the construction of $\G=\mcl_{4m,2}$, it is easy to see that in $\G=\mcl_{4m,2}$ passing each vertex there is exactly one $4$-cycle, which belongs to $\mathcal{F}$. Clearly, any two distinct $4$-cycles in $\mathcal{F}$ are vertex-disjoint. This implies that $\Delta=\{V(\C_{j}^{i})\ |\ i\in\mz_{2},j\in\mz_{m}\}$ is an $A$-invariant partition of $V(\G)$. Consider the quotient graph $\G_{\Delta}$, and let $T$ be the kernel of $A$ acting on $\Delta$. Then $\G_{\Delta}\cong \C_{m}[2K_{1}]$, the lexicographic product of a cycle of length $m$ and an empty graph of order $2$. Hence $A/T\leq \Aut(\C_{m}[2K_{1}])\cong \mz_{2}^{m}\rtimes D_{2m}$. Note that between any two adjacent vertices of $\G_{\Delta}$ there is exactly one edge of $\G=\mcl_{4m,2}$. Then $T$ fixes each vertex of $\G$ and hence $T=1$.
So we may view $A$ as a subgroup of $\Aut(\G_\Delta)\cong\Aut(\C_{m}[2K_{1}])\cong\mz_{2}^{m}\rtimes D_{2m}$.

For convenience, we will simply use the $\C_{j}^{i}$'s to represent the vertices of $\G_{\Delta}$.
Then $\G_{\Delta}$ has vertex set
\[\{\C_{j}^{0}, \C_{j}^{1}\ |\ j\in\mz_{m}\}\]
and edge set
\[\{\{\C_{j}^{0}, \C_{j+1}^{0}\}, \{\C_{j}^{1}, \C_{j+1}^{1}\}, \{\C_{j}^{0}, \C_{j+1}^{1}\}, \{\C_{j}^{1}, \C_{j+1}^{0}\} \ |\ j\in\mz_{m}\}.\]

Let $\mb=\{\{\C_{j}^{0}, \C_{j}^{1}\}\ |\ j\in\mz_{m}\}$. Then $\mb$ is an $\Aut(\G_\Delta)$-invariant partition of $V(\G_{\Delta})$. Let $K$ be the kernel of $\Aut(\G_\Delta)$ acting on $\mb$.
Then $K=\lg k_{0}\rg \times \lg k_{2}\rg\times  \cdots\times \lg k_{m-1}\rg$, where we use $k_{i}$ to denote the transposition $(\C_{j}^{0}\ \C_{j}^{1})$ for $j\in\mz_{m}$. Clearly, $K$ is the maximal normal $2$-subgroup of $\Aut(\G_\Delta)$.

Suppose to the contrary that $\G=\mcl_{4m,2}$ is a Cayley graph. By Proposition~\ref{IsCay}, $A$ has a subgroup, say $G$ acting regularly on $V(\G)$. Then $G$ has order $8m$, and \[G/(G\cap K)\cong GK/K\leq \Aut(\G_\Delta)/K\lesssim D_{2m}.\] Since $m$ odd, it follows that $|G\cap K|=4$ or $8$, and so
$G\cap K\cong \mz_{2}^{2}$ or $\mz_{2}^{3}$.

If $G\cap K\cong \mz_{2}^{2}$, then $|GK/K|=2m$ and $GK/K=  \Aut(\G_\Delta)/K\cong D_{2m}$. So $GK= \Aut(\G_\Delta)\cong\mz_{2}^{m}\rtimes D_{2m}$.
Let $M$ be a Hall $2'$-subgroup of $G$. Then $M\cong\mz_m$ and $M$ is also a Hall $2'$-subgroup of $\Aut(\G_\Delta)$. Clearly, $\Aut(\G_\Delta)$ is solvable, so all Hall $2'$-subgroups of $\Aut(\G_\Delta)$ are conjugate. Without loss of generality, we may
let $M=\lg\a\rg$, where $\a$ is the following permutation on $V(\G_\Delta)$:
\[\a=(\C_{0}^{0}~\C_{1}^{0}\ldots \C_{m-1}^{0})(\C_{0}^{1}~\C_{1}^{1}\ldots \C_{m-1}^{1}).\]
Then $K\rtimes\lg\a\rg$ acts transitively on $V(\G_\Delta)$. Clearly, $C_K(\a)$ is contained in the center of $K\rtimes\lg\a\rg$.
So $C_K(\a)$ is semiregular on $V(\G_\Delta)$. This implies that
\[C_K(\a)=\lg k_{0}k_{1}\ldots k_{m-1}\rg\cong \mz_{2}.\]
On the other hand, let $L=(G\cap K)M$. Clearly, $G\cap K\unlhd G$, so $L$ is a subgroup of $G$ of order $4m$. {For any odd prime factor $p$ of $m$}, let $P$ be a Sylow $p$-subgroup of $M$. Then $P$ is also a Sylow $p$-subgroup of $L$, and since $M$ is cyclic, one has $M\leq N_L(P)$. By Sylow theorem, we have $|L: N_L(P)|=kp+1\mid 4$ for some integer $k$. Since $3\nmid m$, one has $L=N_L(P)$. It follows that $M\unlhd L$ and so $L=M\times (G\cap K)$. This implies that $G\cap K\leq C_{K}(M)=C_{K}(\a)\cong\mz_2$, a contradiction.

If $G\cap K\cong \mz_{2}^{3}$, then $|GK/K|=m$. Furthermore, $GK/K\cong\mz_{m}$ and $GK/K$ acts on $\mb$ regularly.
Since $G$ is transitive on $V(\G)$, there exists $g\in G$ such that $(x_{0}^{1,1})^{g}=x_{1}^{1,1}$, where $x_{0}^{1,1}, x_{1}^{1,1}\in \C_{0}^{1}$.
As $V(\G_\Delta)=\{\C_{j}^{i}\ |\ i\in\mz_{2}, j\in\mz_{m}\}$, $g$ fixes the $4$-cycle $\C_{0}^{1}=(x_{0}^{1,1}, x_{1}^{1,1}, x_{0}^{1,0}, x_{1}^{1,0})$. Since $\mb=\{\{\C_{j}^{0}, \C_{j}^{1}\}\ |\ j\in\mz_{m}\}$ is also $A$-invariant, $g$ fixes $\{\C_{0}^{0}, \C_{0}^{1}\}$ setwise. Since $GK/K$ acts on $\mb$ regularly,  $g$ fixes $\{\C_{j}^{0}, \C_{j}^{1}\}$ setwise for every $j\in\mz_{m}$.
Observe that $\{x_{0}^{1,1}, x_{2m-1}^{1,1}\}$ and $\{x_{1}^{1,1}, x_{2}^{1,1}\}$ are the unique edges of $\G$ between $\C_{0}^{1}$ and $\C_{m-1}^{1}$, $\C_{0}^{1}$ and $\C_{2}^{1}$, respectively. This implies that $g$ will map $\C_{m-1}^{1}$ to $\C_{2}^{1}$, contradicting that $g$ fixes $\{\C_{j}^{0}, \C_{j}^{1}\}$ setwise for every $j\in\mz_{m}$.
\hfill\qed

\section{A family of trivalent VNC bi-dihedrants}

The goal of this section is to prove the following lemma which gives a new family of trivalent vertex-transitive non-Cayley bi-dihedrants.
To be brief, a vertex-transitive non-Cayley graph is sometimes simply called a {\em VNC graph}.


\begin{lem}\label{non-Cayley-2}
 Let $H=\lg a,b \ |\ a^n=b^2=1,a^{b}=a^{-1}\rg$ be a dihedral group, where $n=48\ell$ and $\ell\geq 1$. Then $\G=\BiCay(H, \{b, ba\}, \{ba^{24\ell}, ba^{12\ell-1}\}, \{1\})$ is a VNC dihedrant.
\end{lem}

\f\demo We first define a permutation on $V(\G) $ as follows:
$$
\begin{array}{llll}
g: & (a^{3r})_{0} \mapsto (a^{3r})_{0},   &(a^{3r})_{1}\mapsto (ba^{3r})_{0}, &
     (a^{3r+1})_{0}\mapsto (ba^{3r+1})_{1},\\
   & (a^{3r+1})_{1}\mapsto (a^{24\ell+3r+1})_{1},  &(a^{3r+2})_{i}\mapsto (ba^{12\ell+3r+2})_{i+1},&
     (ba^{3r})_{0}\mapsto (a^{3r})_{1},\\
   & (ba^{3r})_{1}\mapsto (ba^{24\ell+3r})_{1}, &(ba^{3r+1})_{0}\mapsto (ba^{3r+1})_{0}, &(ba^{3r+1})_{1} \mapsto(a^{3r+1})_{0},\\
   & (ba^{3r+2})_{i}\mapsto(a^{-12\ell+3r+2})_{i+1},
\end{array}
$$ where $r\in \mz_{16\ell},~i\in\mz_{2}$.

It is easy to check that $g$ is an involution, and furthermore,
for any $t\in \mz_{16\ell}$,
we have
$$
\begin{array}{l}
\G((a^{3r})_{0})^g
=\{(a^{3r})_1, (ba^{3r})_0, (ba^{3r+1})_0\}=\G((a^{3r})_{0}),\\
\G((a^{3r})_{1})^g
=\{(ba^{3r})_1, (a^{3r})_0, (a^{3r-1})_0\}=\G((ba^{3r})_0),\\
\G((ba^{3r})_1)^g
=\{(ba^{24\ell+3r})_0, (a^{3r})_1, (a^{12\ell+3r+1})_1\}=\G((ba^{24\ell+3r})_1),\\
\G((a^{3r+1})_0)^g
=\{(ba^{3r+1})_0, (a^{24\ell+3r+1})_1, (a^{36\ell+3r+2})_1\}=\G((ba^{3r+1})_1),\\
\G((a^{3r+1})_1)^g
=\{(a^{24\ell+3r+1})_0, (ba^{3r+1})_1, (ba^{36\ell+3r})_1\}=\G((a^{24\ell+3r+1})_1),\\
\G((ba^{3r+1})_0)^g
=\{(ba^{3r+1})_1, (a^{3r+1})_0, (a^{3r})_0\}=\G((ba^{3r+1})_0),\\
\G((a^{3r+2})_0)^g
=\{(ba^{12\ell+3r+2})_0, (a^{36\ell+3r+2})_1, (a^{3r+3})_1\}=\G((ba^{12\ell+3r+2})_1),\\
\G((a^{3r+2})_1)^g
=\{(ba^{12\ell+3r+2})_1, (a^{12\ell+3r+2})_0, (a^{12\ell+3r+1})_0\}=\G((ba^{12\ell+3r+2})_0).\\
\end{array}
$$
This implies that $g$ is an automorphism of $\G$. Observing that $g$ maps $1_{1}$ to $b_{0}$, it follows that $\lg \R(H), g\rg$ is transitive on $V(\G)$, and so $\G$ is a vertex-transitive graph.

Below, we shall first prove the following claim.

\medskip

\f{\bf Claim.} $\Aut(\G)_{1_{0}}=\lg g \rg$.

Let $A=\Aut(\G)$. It is easy to see that $g$ fixes $1_{0}$, and so $g\in A_{1_{0}}$. To prove the Claim, it suffices to prove that $|A_{1_{0}}|=2.$

Note that the neighborhood $\G(1_0)$ of $1_0$ in $\G$ is $=\{1_1, b_0, (ba)_0\}$. By a direct computation, we find that in $\G$ there is a unique $8$-cycle passing through $1_{0}$, $1_{1}$ and $b_{0}$, that is,
$$\C_{0}=(1_{0}, 1_{1}, (ba^{24\ell})_{1}, (ba^{24\ell})_{0}, (a^{24\ell})_{0}, (a^{24\ell})_{1}, b_{1}, b_{0}, 1_{0}).$$
Furthermore, in $\G$ there is no $8$-cycle passing through $1_{0}$ and $(ba)_{0}$. So $A_{1_{0}}$ fixes $(ba)_{0}$.

If $A_{1_{0}}$ also fixes $1_{1}$ and $b_{0}$, then $A_{1_{0}}$ will fix every neighbor of $1_0$, and the connectedness and vertex-transitivity of $\G$ give that $A_{1_{0}}=1$, a contradiction.
Therefore, $A_{1_{0}}$ swaps $1_{1}$ and $b_{0}$, and $(ba)_{0}$ is the unique neighbor of $1_0$ such that $A_{1_{0}}=A_{(ba)_{0}}$.
It follows that $\{1_{0}, (ba)_{0}\}$ is a block of imprimitivity of $A$ acting on $V(\G)$.
Since $\G$ is vertex-transitive, every $v\in V(\G)$ has a unique neighbor, say $u$ such that $A_u=A_v$.
Then the set \[\mb=\{\{u, v\}\in E(\G)\mid A_u=A_v\}\]
forms an $A$-invariant partition of $V(\G)$. Clearly, $\{1_{0}, (ba)_{0}\}\in\mb$.
Similarly, since $\C_0$ is also the unique $8$-cycle of $\G$ passing through $1_{0}$, $1_{1}$ and $b_{0}$, $A_{1_1}$ swaps $1_0$ and $b_0$, and $(ba^{12\ell-1})_{1}$ is the unique neighbor of $1_1$ such that $A_{1_{1}}=A_{(ba^{12\ell-1})_{1}}$. So
$\{1_{1}, (ba^{12\ell-1})_{1}\}\in\mb$. Set
\[\mb_0=\{\{1_{0},(ba)_{0}\}^{\R(h)} \ |\ h\in H \}\ {\rm and}\ \mb_1=\{ \{1_{1}, (ba^{12\ell-1})_{1}\}^{\R(h)} \ |\ h\in H \}.\]
Clearly, $\mb=\mb_0\cup\mb_1$.

Now we consider the quotient graph $\G_{\mb}$ of $\G$ relative to $\mb$. It is easy to see that $\lg \R(a)\rg$ acts semiregularly on $\mb$ with $\mb_0$ and $\mb_1$ as its two orbits. So $\G_{\mb}$ is isomorphic to a bi-Cayley graph over $\lg a\rg$.
Set $B_0=\{1_{0}, (ba)_{0}\}$ and $B_1=\{1_{1}, (ba^{12\ell-1})_{1}\}$.
Then one may see that the neighbors of $B_0$ in $\G_{\mb}$ are: $B_0^{\R(a)}, B_0^{\R(a^{-1})}, B_1, B_1^{\R(a^{-12\ell+2})}$,
and the neighbors of $B_1$ in $\G_{\mb}$ are: $B_1^{\R(a^{12\ell+1})}, B_1^{\R(a^{-12\ell-1})}, B_0, B_0^{\R(a^{12\ell-2})}$.
So
\[\G_{\mb}\cong\G'=\BiCay(\lg a\rg, \{a, a^{-1}\}, \{a^{12\ell+1}, a^{-12\ell-1}\}, \{1, a^{-12\ell+2}\}).\]
Observe that there is one and only one edge of $\G$ between $B_0$ and any one of its neighbors in $\G_\mb$.
Clearly, $A$ acts transitively on $V(\G_\mb)$, so there is one and only one edge of $\G$ between every two adjacent blocks of $\mb$.
It follows that $A$ acts faithfully on $V(\G_\mb)$, and hence we may view $A$ as a subgroup of $\Aut(\G_\mb)$.
Recall that $g\in A_{1_0}=A_{(ba)_0}$. Moreover, $g$ swaps the two neighbors $1_1$ and $b_0$ of $1_0$. Clearly, $1_1\in B_1$ and $b_0\in B_0^{\R(a^{-1})}$, so $g$ swaps the two blocks $B_1$ and $B_0^{\R(a^{-1})}$. Similarly, $g$ swaps the two neighbors $(ba)_1$ and $a_0$ of $(ba)_0$.
Clearly, $(ba)_1\in B_1^{\R(a^{-12\ell+2})}$ and $a_0\in B_0^{\R(a)}$, so
$g$ swaps the two blocks $B_1^{\R(a^{-12\ell+2})}$ and $B_0^{\R(a)}$.
Note that $\R(ab)$ swaps the two vertices in $B_0$. So $\lg g, \R(ab)\rg$ acts transitively on the neighborhood of $B_0$ in $\G_\mb$.
This implies that $A$ acts transitively on the arcs of $\G_\mb$, and so $\G'$ is a tetravalent arc-transitive bi-circulant.
In \cite{KKMW}, a characterization of tetravalent edge-transitive bi-circulants is given. It is easy to see that our graph $\G'$ belongs to Class 1(c) of \cite[Theorem~1.1]{KKMW}. By checking \cite[Theorem~4.1]{KKMW}, we see that the stabilizer $\Aut(\G')_{u}$ of $u\in V(\G')$ has order $4$. This implies that $|A|=4|V(\G_\mb)|=8n$. Consequently, $|A_{1_0}|=2$ and so our claim holds.\medskip

Now we are ready to finish the proof. Suppose to the contrary that $\G$ is a Cayley graph.
By Proposition~\ref{IsCay}, $A$ contains a subgroup, say $J$ acting regularly on $V(\G)$.
By Claim, $J$ has index $2$ in $A$, and since $g\in A_{1_0}$, one has $A=J\rtimes \lg g \rg$.
It is easy to check that $\R(a), \R(b)$ and $g$ satisfy the following relations:
$$(g\R(b))^{4}=\R(a^{24\ell}),~g\R(a^{3})=\R(a^{3})g,~g\R(ba)=\R(ba)g,~g=\R(a)(g\R(b))^{2}\R(a^{12\ell-1}).$$

Suppose that $\R(H)\nleq J$. Then $A=J\R(H)$. Since $|J: \R(H)|=2$, it follows that $|\R(H): J\cap\R(H)|=2$. Thus, $J\cap\R(H)=\lg \R(a)\rg$ or $\lg \R(a^2), \R(b)\rg$. If $\R(H)\cap J=\lg \R(a)\rg$, then we have $\R(b)\notin J$, $\R(a)\in J$, and hence $A=J\cup J\R(b)=J\cup Jg$, implying that $J\R(b)=Jg$. It follows that $g\R(b)\in J$, and then $g=\R(a)(g\R(b))^{2}\R(a^{12\ell-1})\in J$ due to $\R(a)\in J$, a contradiction.
If $\R(H)\cap J=\lg \R(a^2), \R(b)\rg$, then $\R(a)\notin J$, and again we have $A=J\cup J\R(a)=J\cup Jg$, implying that $J\R(a)=Jg$.
So, $\R(a)g, g\R(a^{-1})\in J$. Then \[g=\R(a)g\R(b)g\R(b)\R(a^{12\ell-1})=(\R(a)g)\R(b)(g\R(a^{-1}))\R(ba^{12\ell-2})\in J,\] a contradiction.

Suppose that $\R(H)\leq J$. Then $|J: \R(H)|=2$ and $\R(H)\unlhd J$. Since $J$ is regular on $V(\G)$, by Proposition~\ref{normalizer}, there exists a $\delta_{\a, x, y}\in J$ such that $1_{0}^{\delta_{\a, x, y}}=1_{1}$, where $\a\in\Aut(H)$ and $x, y\in H$. By the definition of $\d_{\a,x,y}$, we have $1_1=1_0^{\d_{\a,x,y}}=(x\cdot 1^\a)_1=x_1$, implying that $x=1$. Furthermore, we have the following relations:
\[R^{\alpha}=x^{-1}Lx, L^{\alpha}=y^{-1}Ry, S^{\alpha}=y^{-1}S^{-1}x,\]
where $R=\{b, ba\}, L=\{ba^{24\ell}, ba^{12\ell-1}\}, S=\{1\}$.
In particular, the last equality implies that $x=y$ due to $S=\{1\}$. So we have $x=y=1$.  From the proof of Claim we know that $B_0=\{1_{0}, (ba)_{0}\}$ and $B_1=\{1_{1}, (ba^{12\ell-1})_{1}\}$ are two blocks of imprimitivity of $A$ acting on $V(\G)$. So we have $((ba)_{0})^{\delta_{\a, 1, 1}}=(ba^{12\ell-1})_{1}$. It follows that $(ba)^\a=ba^{12\ell-1}$, and then from $R^\a=L$ we obtain that $b^\a=ba^{24\ell}$. Consequently, we have $a^\a=a^{36\ell-1}$. One the other hand, we have $\{b, ba\}=R=L^\a=\{b, ba^{24\ell+1}\}$. This forces that $ba=ba^{24\ell+1}$, which is clearly impossible.
\hfill\qed



\section{Two families of trivalent Cayley bi-dihedrants}

In this section, we shall prove two lemmas which will be used the proof of Theorem~\ref{all-non-cay}.

\begin{lem}\label{Cayley}
Let $H=\lg a,b \ |\ a^{12m}=b^2=1, a^{b}=a^{-1}\rg$ be a dihedral group with $m$ odd. Then for each $i\in\mz_{12m}$, $\G=\BiCay(H, \{b,ba^{i}\}, \{ba^{6m},ba^{3m-i}\}, \{1\})$ is a Cayley graph whenever $\lg a^i, a^{3m}\rg=\lg a\rg$.
\end{lem}

\f\demo Let $g$ be a permutation of $V(\G)$ defined as follows:
$$
\begin{array}{lll}
 g:&(a^{6km+3ri})_{j}\mapsto (ba^{6(k+1)m+3ri})_{j+1}, &(ba^{6km+3ri})_{j}\mapsto(a^{6km+3ri})_{j+1},\\
   &(a^{3km+(3r+1)i})_{0}\mapsto (a^{3(k+1)m+(3r+1)i})_{0},&(ba^{3km+(3r+1)i})_{0}\mapsto (a^{3(k+1)m+(3r+1)i})_{1},\\
   &(a^{3km+(3r+1)i})_{1}\mapsto (ba^{3(k+1)m+(3r+1)i})_{0}, & (ba^{3km+(3r+1)i})_{1}\mapsto (ba^{3(k-1)m+(3r+1)i})_{1},\\
   &(a^{3km+(3r+2)i})_{0}\mapsto (ba^{3(k+1)m+(3r+2)i})_{1}, &(ba^{3km+(3r+2)i})_{0}\mapsto (ba^{3(k+1)m+(3r+2)i})_{0},\\
   &(a^{3km+(3r+2)i})_{1}\mapsto (a^{3(k-1)m+(3r+2)i})_{1}, & (ba^{3km+(3r+2)i})_{1}\mapsto (a^{3(k+1)m+(3r+2)i})_{0},\\
\end{array}
$$ where $r\in\mz_{m}$, $k\in\mz_{4}$ and $j\in\mz_{2}$.
It is easy to check that for any $r\in\mz_{m}$, $k\in\mz_{4}$ and $j\in\mz_2$, we have
$$
\begin{array}{lll}
\G((a^{6km+3ri})_{0})^g
&=&\{(a^{6km+3ri})_{1}, (ba^{6km+3ri})_{0}, (ba^{6km+(3r+1)i})_{0}\}^{g}\\
&=&\{(ba^{6(k+1)m+3ri})_{0}, (a^{6km+3ri})_{1}, (a^{3(k+1)m+(3r+1)i})_{1}\}\\
&=&\G((ba^{6(k+1)m+3ri})_{1}),\\

\G((a^{6km+3ri})_{1})^g
&=&\{(a^{6km+3ri})_{0}, (ba^{6(k+1)m+3ri})_{1}, (ba^{(6k+3)m+(3r-1)i})_{1}\}^{g}\\
&=&\{(ba^{6(k+1)m+3ri})_{1}, (a^{6(k+1)m+3ri})_{0}, (a^{(6k+6)m+(3r-1)i})_{0}\}\\
&=&\G((ba^{6(k+1)m+3ri})_{0}),\\

\G((a^{3km+(3r+1)i})_{0})^g
&=&\{(a^{3km+(3r+1)i})_{1}, (ba^{3km+(3r+1)i})_{0}, (ba^{3km+(3r+2)i})_{0}\}^{g}\\
&=&\{(ba^{3(k+1)m+(3r+1)i})_{0}, (a^{3(k+1)m+(3r+1)i})_{1}, (ba^{3(k+1)m+(3r+2)i})_{0}\}\\
&=&\G((a^{3(k+1)m+(3r+1)i})_{0}),\\

\G((a^{3km+(3r+1)i})_{1})^g
&=&\{(a^{3km+(3r+1)i})_{0}, (ba^{3(k+2)m+(3r+1)i})_{1}, (ba^{3(k+1)m+3ri})_{1}\}^{g}\\
&=&\{(a^{3(k+1)m+(3r+1)i})_{0}, (ba^{3(k+1)m+(3r+1)i})_{1}, (a^{3(k+1)m+3ri})_{0}\}\\
&=&\G((ba^{3(k+1)m+(3r+1)i})_{0}),\\

\G((a^{3km+(3r+2)i})_{0})^g
&=&\{(a^{3km+(3r+2)i})_{1}, (ba^{3km+(3r+2)i})_{0}, (ba^{3km+3(r+1)i})_{0}\}^{g}\\
&=&\{(a^{3(k-1)m+(3r+2)i})_{1}, (ba^{3(k+1)m+(3r+2)i})_{0}, (a^{3km+3(r+1)i})_{1}\}\\
&=&\G((ba^{3(k+1)m+(3r+2)i})_{1}),\\

\end{array}
$$

$$
\begin{array}{lll}
\G((a^{3km+(3r+2)i})_{1})^g
&=&\{(a^{3km+(3r+2)i})_{0}, (ba^{3(k+2)m+(3r+2)i})_{1}, (ba^{3(k+1)m+(3r+1)i})_{1}\}^{g}\\
&=&\{(ba^{3(k+1)m+(3r+2)i})_{1}, (a^{3(k+3)m+(3r+2)i})_{0}, (ba^{3km+(3r+1)i})_{1}\}\\
&=&\G((a^{3(k-1)m+(3r+2)i})_{1}),\\

\G((ba^{6km+3ri})_{0})^g
&=&\{(ba^{6km+3ri})_{1}, (a^{6km+3ri})_{0}, (a^{6km+(3r-1)i})_{0}\}^{g}\\
&=&\{(a^{6km+3ri})_{0}, (ba^{6(k+1)m+3ri})_{1}, (a^{(6k+3)m+(3r-1)i})_{1}\}\\
&=&\G((a^{6km+3ri})_{1}),\\

\G((ba^{6km+3ri})_{1})^g
&=&\{(ba^{6km+3ri})_{0}, (a^{6(k-1)m+3ri})_{1}, (a^{(6k-3)m+(3r+1)i})_{1}\}^{g}\\
&=&\{(a^{6km+3ri})_{1}, (ba^{6km+3ri})_{0}, (ba^{6km+(3r+1)i})_{0}\}\\
&=&\G((a^{6km+3ri})_{0}),\\

\G((ba^{3km+(3r+1)i})_{0})^g
&=&\{(ba^{3km+(3r+1)i})_{1}, (a^{3km+(3r+1)i})_{0}, (a^{3km+3ri})_{0}\}^{g}\\
&=&\{(ba^{3(k-1)m+(3r+1)i})_{1}, (a^{3(k+1)m+(3r+1)i})_{0}, (ba^{3(k+2)m+3ri})_{1}\}\\
&=&\G((a^{3(k+1)m+(3r+1)i})_{1}),\\

\G((ba^{3km+(3r+1)i})_{1})^g
&=&\{(ba^{3km+(3r+1)i})_{0}, (a^{3(k+2)m+(3r+1)i})_{1}, (a^{3(k-1)m+(3r+2)i})_{1}\}^{g}\\
&=&\{(a^{3(k+1)m+(3r+1)i})_{1}, (ba^{3(k+3)m+(3r+1)i})_{0}, (a^{3(k-2)m+(3r+2)i})_{1}\}\\
&=&\G((ba^{3(k-1)m+(3r+1)i})_{1}),\\

\G((ba^{3km+(3r+2)i})_{0})^g
&=&\{(ba^{3km+(3r+2)i})_{1}, (a^{3km+(3r+2)i})_{0}, (a^{3km+(3r+1)i})_{0}\}^{g}\\
&=&\{(a^{3(k+1)m+(3r+2)i})_{0}, (ba^{3(k+1)m+(3r+2)i})_{1}, (a^{(3k+1)m+3(r+1)i})_{0}\}\\
&=&\G((ba^{3(k+1)m+(3r+2)i})_{0}),\\

\G((ba^{3km+(3r+2)i})_{1})^g
&=&\{(ba^{3km+(3r+2)i})_{0}, (a^{3(k+2)m+(3r+2)i})_{1}, (a^{3(k-1)m+(3r+3)i})_{1}\}^{g}\\
&=&\{(ba^{3(k+1)m+(3r+2)i})_{0}, (a^{3(k+1)m+(3r+2)i})_{1}, (ba^{3(k+1)m+(3r+3)i})_{0}\}\\
&=&\G((a^{3(k+1)m+(3r+2)i})_{0}).
\end{array}
$$
It follows  that $g\in\Aut(\G)$. Furthermore, one may check that $g$ and $\R(a^{2})$ satisfy the following relations:
$$\R(a^{12m})=g^{4}=1,~g^{2}=\R(a^{6m}),~\R(a^{6})g=g\R(a^{6}),~\R(a^{2})g=g\R(a^{4})g\R(a^{-2}).$$
By the last equality, we have \[(\R(a^{2})g)^{3}=[g\R(a^{4})g\R(a^{-2}))]\R(a^{2})g\R(a^{2})g=g\R(a^{4})g^{2}\R(a^{2})g.\]
It then follows from the second and third equalities that
\[g\R(a^{4})g^{2}\R(a^{2})g=g\R(a^{6+6m})g=g^{2}\R(a^{6+6m})=\R(a^{6}).\]
Therefore, $(\R(a^{2})g)^{3}=\R(a^{6}).$

Let $G=\lg \R(a^{2}),\R(b),g\rg$ and $T=\lg \R(a^{6})\rg$. Then $T\unlhd G$ and
\[G/T=\lg \R(a^{2})T, gT \ |\ \R(a^{2})^{3}T=g^{2}T=(\R(a^{2})g)^{3}T=T\rg\cong A_{4}.\] So $|G|=24m$.

Let
$$\begin{array}{ll}
\Omega_{00}=\{t_{0}\ |\ t\in\lg a^{2},b\rg\}, &\Omega_{01}=\{t_{1}\ |\ t\in\lg a^{2},b\rg\},\\ \Omega_{10}=\{(at)_{0}\ |\ t\in\lg a^{2},b\rg\},&\Omega_{11}=\{(at)_{1}\ |\ t\in\lg a^{2},b\rg\}.
\end{array}
$$
Then $\Omega_{ij}$'s $(0\leq i,j\leq 1)$ are orbits of $T$ and $V(\G)=\bigcup\limits_{0\leq i,j\leq1}\Omega_{ij}$. Since $1_{0}^{g}=(ba^{6m})_{1}\in \Omega_{01}$, $a_{0}^{g}=(a^{3m+1})_{0}\in \Omega_{00}$ and $a_{1}^{g}=(ba^{3m+1})_{1}\in\Omega_{01}$, it follows that $G$ is transitive, and so regular on $V(\G)$. By Proposition~\ref{IsCay}, $\G$ is a Cayley graph on $G$, as required.
\hfill\qed

\begin{lem}\label{Cayley2}
Let $H=\lg a,b \ |\ a^{12m}=b^2=1,a^{b}=a^{-1}\rg$ be a dihedral group with $m$ even and $4\nmid m$. Then the following two bi-Cayley graphs:
\[\begin{array}{l}
\G_{1}=\BiCay(H, \{b,ba\}, \{ba^{6m}, ba^{3m-1}\}, \{1\}),\\
 \G_{2}=\BiCay(H, \{b,ba\}, \{ba^{6m}, ba^{9m-1}\},\{1\})
 \end{array}\]
are both Cayley graphs.
\end{lem}

\f\demo Let $V=H_0\cup H_1$. Then $V(\G_1)=V(\G_2)=V$. We first define two permutations on $V$ as follows:
$$
\begin{array}{lll}
 g_{1}: &(a^{4r})_{i}\mapsto(ba^{6m+4r})_{i+1}, &(ba^{4r})_{i}\mapsto (a^{4r})_{i+1},\\
        &(a^{4r+1})_{i}\mapsto (ba^{9m+4r+1})_{i+1}, &(ba^{4r+1})_{i}\mapsto (a^{3m+4r+1})_{i+1},\\
        &(a^{4r+2})_{i}\mapsto (ba^{4r+2})_{i+1}, &(ba^{4r+2})_{i}\mapsto (a^{6m+4r+2})_{i+1},\\
        &(a^{4r+3})_{i}\mapsto (ba^{3m+4r+3})_{i+1}, &(ba^{4r+3})_{i}\mapsto (a^{9m+4r+3})_{i+1},
\end{array}
$$
$$
\begin{array}{lll}
 g_{2}: &(a^{4r})_{i}\mapsto(ba^{6m+4r})_{i+1}, &(ba^{4r})_{i}\mapsto (a^{4r})_{i+1},\\
        &(a^{4r+1})_{i}\mapsto (ba^{3m+4r+1})_{i+1}, &(ba^{4r+1})_{i}\mapsto (a^{9m+4r+1})_{i+1},\\
        &(a^{4r+2})_{i}\mapsto (ba^{4r+2})_{i+1}, &(ba^{4r+2})_{i}\mapsto (a^{6m+4r+2})_{i+1},\\
        &(a^{4r+3})_{i}\mapsto (ba^{9m+4r+3})_{i+1}, &(ba^{4r+3})_{i}\mapsto (a^{3m+4r+3})_{i+1},
\end{array}
$$ where $r\in\mz_{3m}$ and $i\in\mz_{2}$.

It is easy to check that $g_{j}\in\Aut(\G_{j})$ for $j=1$ or $2$. Furthermore, $\R(a^{2}), \R(b) $ and $ g_{j}$ $(j=1$ or $2)$ satisfy the following relations:
$$\begin{array}{l}
\R(a^{12m})=\R(b^{2})=g_{j}^{4}=1, \R(b)\R(a^{2})\R(b)=\R(a^{-2}),\\
 g_{j}^{2}=\R(a^{6m}), \R(b)g_{j}{\R(b)}=g_{j}^{-1}, \\
g_1^{-1}\R(a){g_{1}}=\R(a^{3m+1}), g_2^{-1}\R(a){g_{2}}=\R(a^{9m+1}).
\end{array}$$

For $j=1$ or $2$, let $G_{j}=\lg \R(a),\R(b),g_{j}\rg$.
From the above relations it is east to see that
\[G_{j}=(\lg \R(a)\rg\lg g_{j}\rg)\rtimes\lg \R(b)\rg\]
has order at most $48m$. Observe that $1_{0}^{g_{j}}=(ba^{6m})_{1}\in H_{1}$ for $j=1$ or $2$. It follows that $G_{j}$ is transitive on $V(\G_{j})$, and so $G_{j}$ acts regularly on $V(\G_{j})$. By Proposition~\ref{IsCay}, each $\G_{j}$ is a Cayley graph.
\hfill\qed

\section{Vertex-transitive trivalent bi-dihedrants}

In this section, we shall give a complete classification of trivalent vertex-transitive non-Cayley bi-dihedrants.
For convenience of the statement, throughout this section, we shall make the following assumption.\medskip

\f{\bf Assumption~I.}
\begin{itemize}
  \item $H$: the dihedral group $D_{2n}=\lg a, b\ |\ a^n=b^2=1, bab=a^{-1}\rg (n\geq 3)$,
  \item $\G={\rm BiCay}(H,~R,~L,~\{1\})$: a connected trivalent $2$-type vertex-transitive bi-Cayley graph over the group $H$ (in this case, $|R|=|L|=2$),
  \item $G$: a minimum group of automorphisms of $\G$ subject to that $\R(H)\leq G$ and $G$ is transitive on the vertices but intransitive on the  arcs of $\G$.
\end{itemize}

The following lemma given in \cite{ZZ} shows that the group $G$ must be solvable.

\begin{lem}{\rm\cite[Lemma 6.2]{ZZ}}\label{solvable}
$G=\R(H)P$ is solvable, where $P$ is a Sylow $2$-subgroup of $G$.
\end{lem}

\subsection{$H_{0}$ and $H_{1}$ are blocks of imprimitivity of $G$}

The case where $H_{0}$ and $H_{1}$ are blocks of imprimitivity of $G$ has been considered in \cite{ZZ}, and the main result is the following proposition.

\begin{prop}{\rm\cite[Theorem 1.3]{ZZ}}\label{block}
If $H_0$ and $H_1$ are blocks of imprimitivity of $G$ on $V(\G)$,
then either $\G$ is Cayley or one of the following occurs:
\begin{enumerate}
\item [$(1)$]  $(R, L, S)\equiv(\{b,~ba^{\ell+1}\}, \{ba,~ba^{\ell^{2}+\ell+1}\}, \{1\})$, where $n\geq 5$,
$\ell^{3}+\ell^{2}+\ell+1\equiv0~(\mod n)$, $\ell^{2}\not\equiv 1~(\mod n)$;
\item [$(2)$]  $(R, L, S)\equiv(\{ba^{-\ell},~ba^{\ell}\}, \{a,~a^{-1}\}, \{1\})$, where $n=2k$ and $\ell^{2}\equiv-1~(\mod k)$. Furthermore, $\G$ is also a bi-Cayley graph over an abelian group $\mathbb{Z}_{n}\times \mathbb{Z}_{2}$.
\end{enumerate}
Furthermore, all of the graphs arising from (1)-(2) are vertex-transitive non-Cayley.
\end{prop}

In particular, it is proved in \cite{ZZ} that if $n$ is odd and $\G$ is not a Cayley graph, then $H_{0}$ and $H_{1}$ are blocks of imprimitivity of $G$ on $V(\Gamma)$. Consequently, we can get a classification of trivalent vertex-transitive non-Cayley bi-Cayley graphs over a dihedral group $D_{2n}$ with $n$ odd.

\begin{prop}{\rm\cite[Proposition 6.4]{ZZ}}\label{n-is-odd}
If $n$ is odd, then either $\G$ is a Cayley graph, or $H_{0}$ and $H_{1}$ are blocks of imprimitivity of $G$ on $V(\Gamma)$.
\end{prop}

\subsection{$H_{0}$ and $H_{1}$ are not blocks of imprimitivity of $G$}

In this subsection, we shall consider the case where $H_{0}$ and $H_{1}$ are not blocks of imprimitivity of $G$ on $V(\G)$. We begin by citing a lemma from \cite{ZZ}.

\begin{lem}{\rm\cite[Lemma 6.3]{ZZ}}\label{quotient}
Suppose that $H_{0}$ and $H_{1}$ are not blocks of imprimitivity of $G$ on $V(\G)$. Let $N$ be a normal subgroup of $G$, and let $K$ be the kernel of $G$ acting on $V(\G_{N})$. Let $\Delta$ be an orbit of $N$. If $N$ fixes $H_{0}$ setwise, then one of the following holds:
\begin{enumerate}
\item [$(1)$] $\G[\Delta]$ has valency $1$, $|V(\G_{N})|\geq 3$ and $\G$ is a Cayley graph;

\item [$(2)$] $\G[\Delta]$ has valency $0$, $\G_{N}$ has valency $3$, and $K=N$ is semiregular.
\end{enumerate}
\end{lem}

The following lemma deals with the case where $\Core_G(\R(H))=1$, and in this case we shall see that $\G$ is just the cross ladder graph.

\begin{lem}\label{core}
Suppose that $H_{0}$ and $H_{1}$ are not blocks of imprimitivity of $G$ on $V(\G)$. If $\Core_{G}(\R(H))=\bigcap _{g\in G}\R(H)=1$,
then $\G$ is isomorphic to the cross ladder graph $\cl_{4n}$ with $n$ odd, and furthermore, for any minimal normal subgroup $N$ of $G$, we have the following:
\begin{enumerate}
  \item [{\rm (1)}]\ $N$ is a $2$-group which is non-regular on $V(\G)$;
  \item [{\rm (2)}]\ $N$ does not fix $H_{0}$ setwise;
  \item [{\rm (3)}]\ every orbit of $N$ consists of two non-adjacent vertices.
\end{enumerate}
\end{lem}

\f\demo Let $N$ be a minimal normal subgroup of $G$. By Lemma~\ref{solvable}, $G$ is solvable. It follows that $N$ is an elementary abelian $r$-subgroup for some prime divisor $r$ of $|G|$. Clearly, $N\nleq \R(H)$ due to $\Core_{G}(\R(H))=1$. Then $|N\R(H)|/|\R(H)|\mid |G|/|\R(H)|$. From Lemma~\ref{solvable} it follows that $|G|/|\R(H)|$ is a power of $2$, and hence $N$ is a $2$-group.

{Suppose that $N$ is regular on $V(\G)$. Then $N\R(H)$ is transitive on $V(\G)$ and $\R(H)$ is also a $2$-group. Therefore, $N\R(H)$ is not transitive on the arcs of $\G$. The minimality of $G$ gives that $G=N\R(H)$.
Since $n$ is even, $\R(a^{\frac{n}{2}})$ is in the center of $\R(H)$. Set $Q=N\lg \R(a^{\frac{n}{2}})\rg$. Then $Q\unlhd G$ and then $1\neq N\cap Z(Q)\unlhd G$. Since $N$ is a minimal normal subgroup of $G$, one has $N\leq Z(Q)$, and hence $Q$ is abelian. It follows that $\lg \R(a^{\frac{n}{2}})\rg\unlhd G$, contrary to the assumption that $\Core_{G}(\R(H))=1$. Thus, $N$ is not regular on $V(\G)$.
(1) is proved.}

For (2), by way of contradiction, suppose that $N$ fixes $H_{0}$ setwise. Consider the quotient graph $\G_{N}$ of $\G$ relative to $N$, and
let $K$ be the kernel of $G$ acting on $V(\G_{N})$. Take $\Delta$ to be an orbit of $N$ on $V(\G)$.
Then either $(1)$ or $(2)$ of Lemma~\ref{quotient} happens.

For the former, $\G[\Delta]$ has valency $1$ and $|V(\G_{N})|\geq 3$. Then $\G_N$ is a cycle. Moreover, any two neighbors of $u\in\Delta$ are in different orbits of $N$. It follows that the stabilizer $N_{v}$ of $v$ in $N$ fixes every neighbor of $u$. The connectedness of $\G$ implies that $N_{v}=1$. Thus, $K=N$ is semiregular and $\G_{N}$ is a cycle of length $\ell=2|\R(H)|/|N|$. So $G/N\leq \Aut(\G_{N})\cong D_{2\ell}$. If $G/N<\Aut(\G_{N})$, then $|G:N|=\ell$ and so $|G|=2|\R(H)|$. This implies that $\R(H)\unlhd G$, contrary to the assumption that $\Core_{G}(\R(H))=1$. If $G/N=\Aut(\G_{N})$, then $|G: \R(H)|=4$. Since $N\not\leq \R(H)$ and since $N$ fixes $H_{0}$ setwise,
one has $|G: \R(H)N|=2$. It follows that $\R(H)N\unlhd G$. Clearly, $H_{0}$ and $H_{1}$ are just two orbits of $\R(H)N$, and they are also two blocks of imprimitivity of $G$ on $V(\G)$, a contradiction.

For the latter, $\G[\Delta]$ has valency $0$, $\G_{N}$ has valency $3$ and $N=K$ is semiregular.
Let $\bar{H}_{i}$ be the set of orbits of $N$ contained in $H_{i}$ with $i=1,2$. Then $\G_{N}[\bar{H}_{0}]$ and $\G_{N}[\bar{H}_{1}]$ are of valency $2$ and the edges between $\bar{H}_{0}$ and $\bar{H}_{1}$ form a perfect matching. Without loss of generality, we may assume that $1_{0}\in \Delta$. Since $\R(H)$ acts on $H_0$ by right multiplication, we have the subgroup of $\R(H)$ fixing $\Delta$ setwise is just $\R(H)_{\Delta}=\{\R(h)\ |\ h_{0}\in\Delta\}$. If $\R(H)_{\Delta}\leq \lg \R(a)\rg$, then $\R(H)_{\Delta}\unlhd\R(H)$, and the transitivity of $\R(H)$ on $H_0$ implies that $\R(H)_{\Delta}$ will fix all orbits of $N$ contained in $H_{0}$. Since the edges between $\bar{H}_{0}$ and $\bar{H}_{1}$ are independent, $\R(H)_{\Delta}$ fixes all orbits of $N$. It follows that $\R(H)_{\Delta}\leq N$, namely, $\R(H)N/N$ acts regularly on $\bar{H}_{0}$. Then $|\R(H)/(\R(H)\cap N)|=|\R(H)N/N|=|H_{0}/N|$, and so $|N|=|\R(H)\cap N|$, forcing $N\leq \R(H)$, a contradiction. Thus, $\R(H)_{\Delta}\nleq \lg \R(a)\rg$, and so $\lg \R(a)\rg\R(H)_{\Delta}=\R(H)$. This implies that $\lg \R(a), N\rg/N$ is transitive and so regular on $\bar{H}_{0}$. Similarly, $\lg \R(a), N\rg/N$ is also regular on $\bar{H}_{1}$. Thus, $\G_{N}$ is a trivalent $2$-type bi-Cayley graph over $\lg\R(a), N\rg/N$. By \cite[Lemma~5.3]{ZF}, $\bar{H}_{0}$ and $\bar{H}_{0}$ are blocks of imprimitivity of $G/N$, and so $H_{0}$ and $H_{1}$ are blocks of imprimitivity of $G$, a contradiction.

So far, we have completed the proof of (2). Then $N$ does not fix $H_{0}$ setwise, and then $N\R(H)$ is transitive on $V(\G)$. The minimality of $G$ gives that $G=N\R(H)$. Let $P$ and $P_{1}$ be Sylow $2$-subgroups of $G$ and $\R(H)$, respectively, such that $P_{1}\leq P$. Then $N\leq P$ and $P=NP_{1}$.

{If $n$ is even, then by a similar argument to the second paragraph, a contradiction occurs}. Thus, $n$ is odd.
As $H\cong D_{2n}$, $P_{1}\cong\mz_{2}$ and $P_1$ is non-normal in $\R(H)$. So $N\cap \R(H)=1$. Clearly, $|V(\G)|=4n$. If $N$ is semiregular on $V(\G)$, then $N\cong \mz_{2}$ or $\mz_{2}\times\mz_{2}$, and then $|G|=|\R(H)||N|=2|\R(H)|$ or $4|\R(H)|$. Since $\Core_G(\R(H))=1$, we must have
$|G: \R(H)|=4$ and $G\lesssim \Sym(4)$. Since $n$ is odd, one has $n=3$ and $H\cong \Sym(3)$. So $G\cong \Sym(4)$ and hence $G_{1_{0}}\cong \mz_{2}$. Then all involutions of $G (\cong\Sym(4))$ not contained in $N$ are conjugate. Take $1\neq g\in G_{1_{0}}$. Then $g$ is an involution which is not contained in $N$ because $N$ is semiregular on $V(\G)$. Since $\R(H)\cap N=1$, every involution in $\R(H)$ would be conjugate to $g$. This is clearly impossible because $\R(H)$ is semiregular on $V(\G)$. Thus, $N$ is not semiregular on $V(\G)$. (3) is proved.

{Since $n$ is odd, we have $|V(\G_{N})|>2$.} Since $N$ is not semiregular on $V(\G)$, $\G_{N}$ has valency $2$ and $\G[\Delta]$ has valency $0$. This implies that the subgraph induced by any two adjacent two orbits of $N$ is either a union of several cycles or a perfect matching. Thus, $\G_{N}$ has even order. As $\G$ has order $4n$ with $n$ odd, every orbit of $N$ has length $2$. It is easy to see that $\G$ is isomorphic to the cross ladder graph $\cl_{4n}$.
\hfill\qed

The following is the main result of this section.

\begin{theorem}\label{not-block}
Suppose that $H_{0}$ and $H_{1}$ are not blocks of imprimitivity of $G$ on $V(\G)$. Then $\G=\BiCay(H, R,L,S)$ is vertex-transitive non-Cayley if and only if one of the followings occurs:
\begin{enumerate}
\item [$(1)$] $(R,L,S)\equiv (\{b, ba\}, \{b, ba^{2m}\}, \{1\})$, where $n=2(2m+1)$, $m\not\equiv 1~(\mod 3)$, and the corresponding graph is isomorphic the multi-cross ladder graph $\mcl_{4m,2}$;

\item [$(2)$] $(R,L,S)\equiv (\{b, ba\}, \{ba^{24\ell}, ba^{12\ell-1}\}, \{1\})$, where $n=48\ell$ and $\ell\geq 1$.
\end{enumerate}
\end{theorem}

\f\demo The sufficiency can be obtained from Theorem~\ref{EX(m,2)-Cayley} and Lemma~\ref{non-Cayley-2}. We shall prove the necessity in the following subsection by a series of lemmas.\hfill\qed

\subsection{Proof of the necessity of Theorem~\ref{not-block}}

The purpose of this subsection is to prove the necessity of Theorem~\ref{not-block}. Throughout this subsection, we shall always assume that $H_{0}$ and $H_{1}$ are not blocks of imprimitivity of $G$ on $V(\G)$ and that $\G=\BiCay(H, R, L, S)$ is vertex-transitive non-Cayley. In this subsection, we shall always use the following notation.

\medskip
\f{\bf Assumption~II.}\ Let $N=\Core_{G}(\R(H))$. \medskip

Our first lemma gives some properties of the group $N$.

\begin{lem}\label{lem-not-block-1}
$1<N<\lg\R(a)\rg$, $|\lg\R(a)\rg:  N|=n/|N|$ is odd and the quotient graph $\G_N$ of $\G$ relative to $N$ is isomorphic to the cross ladder graph $\cl_{4n/|N|}$.
\end{lem}

\f\demo\  If $N=1$, then from Lemma~\ref{core} it follows that $\G\cong\cl_{4n}$ which is a Cayley graph by Theorem~\ref{not-at-dihedrant}, a contradiction. Thus, $N>1$.
Since $H_{0}$ and $H_{1}$ are not blocks of imprimitivity of $G$ on $V(\G)$, one has $N<\R(H)$.

Consider the quotient graph $\G_{N}$. Clearly, $N$ fixes $H_{0}$ setwise. Recall that $H_{0}$ and $H_{1}$ are not blocks of imprimitivity of $G$ on $V(\G)$ and that $\G$ is non-Cayley. Applying Lemma~\ref{quotient}, we see that $\G_{N}$ is a trivalent $2$-type bi-Cayley graph over $\R(H)/N$. This implies that $|\R(H): N|>2$, and since $H$ is a dihedral group, one has $N<\lg \R(a)\rg$.

Again, by Lemma~\ref{quotient}, $\R(H)/N$ acts semiregularly on $V(\G_{N})$ with two orbits, $\bar{H}_{0}$ and $\bar{H}_{1}$, where $\bar{H}_{i}$ is the set of orbits of $N$ contained in $H_{i}$ with $i=1,0$. Furthermore, $N$ is just the kernel of $G$ acting on $V(\G_N)$ and $N$ acts semiregularly on $V(\G)$. Then $G/N$ is also a minimal vertex-transitive automorphism group of $\G_{N}$ containing $\R(H)/N$. If $\bar{H}_{0}$ and $\bar{H}_{1}$ are blocks of imprimitivity of $G/N$ on $V(\G_N)$, then $H_{0}$ and $H_{1}$ will be blocks of imprimitivity of $G$ on $V(\G)$, which is impossible by our assumption. Thus, $\bar{H}_{0}$ and $\bar{H}_{1}$ are not blocks of imprimitivity of $G/N$ on $V(\G_N)$. Since $N=\Core_{G}(\R(H))$, $\Core_{G/N}(\R(H)/N)$ is trivial. Then from Lemma~\ref{core} it follows that $\G_{N}\cong \cl_{\frac{4n}{|N|}}$, where $\frac{n}{|N|}$ is odd.\hfill\qed

Next, we introduce another notation which will be used in the proof.

\medskip
\f{\bf Assumption~III.}\ Take $M/N$ to be a minimal normal subgroup of $G/N$. \medskip

We shall first consider some basic properties of the quotient graph $\G_M$ of $\G$ relative to $M$.

\begin{lem}\label{lem-not-block-2}
The quotient graph $\G_M$ of $\G$ relative to $M$ is a cycle of length $n/|N|$.
Furthermore, every orbit of $M$ on $V(\G)$ is a union of an orbit of $N$ on $H_{0}$ and an orbit of $N$ on $H_{1}$, and these two orbits of $N$ are non-adjacent.
\end{lem}

\f\demo\ Applying Lemma~\ref{core} to $\G_N$ and $G/N$, we obtain the following facts:
\begin{enumerate}
  \item [{\rm (a)}]\ $M/N$ is an elementary abelian $2$-group which is not regular on $V(\G_N)$,
  \item [{\rm (b)}]\ $M/N$ does not fix $\bar{H}_{0}$ setwise,
  \item [{\rm (c)}]\ every orbit of $M/N$ on $V(\G_N)$ consists of two non-adjacent vertices of $\G_N$.
\end{enumerate}
From (b) and (c) it follows that every orbit of $M$ on $V(\G)$ is just a union of an orbit of $N$ on $H_{0}$ and an orbit of $N$ on $H_{1}$, and these two orbits are non-adjacent. Since every orbit of $N$ on $V(\G)$ is an independent subset of $V(\G)$, each orbit of $M$ on $V(\G)$ is also an independent subset.

Recall that $\G_{N}\cong \cl_{4m}$ where $m=\frac{n}{|N|}$ is odd. The quotient graph of $\G_N$ relative to $M/N$ is just a cycle of length $m$, and so the quotient graph $\G_M$ of $\G$ relative to $M$ is also a cycle of length $m$. \hfill\qed

By Lemma~\ref{lem-not-block-2}, each orbit of $M$ on $V(\G)$ is an independent subset. It follows that the subgraph induced by any two adjacent orbits of $M$ is either a perfect matching or a union of several cycles. For convenience of the statement, the following notations will be used in the remainder of the proof:\medskip

\f{\bf Assumption~IV.}\ \begin{enumerate}
                             \item [{\rm (1)}]\ Let $\Delta$ and $\Delta'$ be two adjacent orbits of $M$ on $V(\G)$ such that $\G[\Delta\cup\Delta']$ is a  union of several cycles.
                             \item [{\rm (2)}]\ Let $\Delta=\Delta_{0}\cup\Delta_{1}$ and $\Delta'=\Delta_{0}'\cup\Delta_{1}'$, where $\Delta_{0},\Delta_{0}'\subseteq H_{0}$ and $\Delta_{1},\Delta_{1}'\subseteq H_{1}$ are four orbits of $N$ on $V(\G)$.
                                 \item [{\rm (3)}]\ $1_{0}\in\Delta_{0}$.
                           \end{enumerate}

Since $\G[\Delta]$ and $\G[\Delta']$ are both null graphs and since $\G[\Delta\cup\Delta']$ is a union of several cycles,  we have the following easy observation.

\begin{lem}\label{lem-not-block-3}
$\G[\Delta_{i}\cup\Delta'_{j}]$ is a perfect matching for any $0\leq i,j\leq 1$.
\end{lem}

The following lemma tells us the possibility of $R$ (Recall that we assume that $\G=\BiCay(H, R, L, \{1\})$).

\begin{lem}\label{lem-not-block-4}
Up to graph isomorphism, we may assume that $R=\{b, ba^i\}$ with $i\in\mz_n\setminus\{0\}$ and that $b_0\in \Delta_{0}'$. Furthermore, we have
\[\begin{array}{l}
\Delta_{0}=\{h_{0}\ |\ \R(h)\in N\}, \Delta'_{0}=\{(bh)_{0}\ |\ \R(h)\in N\}, \\
\Delta'_{1}=\{h_{1}\ |\ \R(h)\in N\},  \Delta_{1}=\{(bh)_{1}\ |\ \R(h)\in N\},
\end{array}
\]
and $1_{1}$ is adjacent to $(ba^{l})_{1}\in \Delta_{1}$ for some $\R(a^l)\in N$.
\end{lem}

\f\demo\ Recall that $N$ is a proper subgroup of $\lg\R(a)\rg$ and that $n/|N|$ is odd. Since $n$ is even by Proposition~\ref{n-is-odd}, it follows that $N$ is of even order, and so the unique involution $\R(a^{n/2})$ of $\lg\R(a)\rg$ is contained in $N$. As $1_0\in \Delta_{0}$ and $N\leq\lg\R(a)\rg$ acts on $H_0$ by right multiplication, one has $\Delta_{0}=\{h_0\mid h\in N\}$. Since $\G[\Delta_{0}]$ is an empty graph, one has $a^{n/2}\notin R$. By Proposition~\ref{properties}~(1), we have $\lg R\cup L\rg=H$, and since $R$ and $L$ are both self-inverse, either $R\subseteq b\lg a\rg$ or $L\subseteq b\lg a\rg$. By Proposition~\ref{properties}~(4), we may assume that $R\subseteq b\lg a\rg$.

Recall that $\G[\Delta_{i}\cup\Delta'_{j}]$ is a perfect matching for any $0\leq i,j\leq 1$. Then $1_0$ is adjacent to $r_0\in \Delta'_{0}$ for some $r\in R$. Since $R\subseteq b\lg a\rg$ and $\Aut(H)$ is transitive on $b\lg a\rg$, by Proposition~\ref{properties}~(3), we may assume that $r=b.$ So $1_{0}$ is adjacent to $b_{0}\in \Delta'_{0}$. Since $N\leq\lg\R(a)\rg$ acts on $H_i$ with $i=0$ or $1$ by right multiplication, we see that the two orbits $\Delta_{0}, \Delta_{0}'$ of $N$ are just the form as given in the lemma. Since $S=\{1\}$, the edges between $H_0$ and $H_1$ form a perfect matching. This enables us to obtain another two orbits $\Delta_{1}, \Delta_{1}'$ of $N$ which have the form as given in the lemma.

By Lemma~\ref{lem-not-block-3}, $\G[\Delta_{1}\cup\Delta'_{1}]$ is a perfect matching. So we may assume that $1_{1}$ is adjacent to $(ba^{l})_{1}\in \Delta_{1}$ for some $\R(a^l)\in N$. \hfill\qed

Now we shall introduce some new notations which will be used in the following.\medskip

\f{\bf Assumption~V.}\ \begin{enumerate}
                             \item [{\rm (1)}]\ Let $T=\lg \R(a^{l})\rg$ be of order $t$, where $a^l$ is given in the above lemma.
                             \item [{\rm (2)}]\ Let
                             $$\begin{array}{ll}
   \Omega_{0}=\{(a^{i\frac{n}{t}})_{0}\ |\ 0\leq i\leq t-1\},&\Omega_{1}=\{(ba^{i\frac{n}{t}})_{1}\ |\ 0\leq i\leq t-1\}, \\
   \Omega'_{0}=\{(ba^{i\frac{n}{t}})_{0}\ |\ 0\leq i\leq t-1\},&\Omega'_{1}=\{(a^{i\frac{n}{t}})_{1}\ |\ 0\leq i\leq t-1\}.
\end{array}
$$
                                 \item [{\rm (3)}]\ $\mb=\{B^{\R(h)}\ |\ h\in H\}$, where $B=\Omega_{0}\cup\Omega_{1}$.
                                 \item [{\rm (4)}]\  Let $B'=\Omega'_{0}\cup\Omega'_{1}$. Then $B'=B^{\R(b)}$.
                           \end{enumerate}

\begin{lem}\label{lem-not-block-5}
The followings hold.
\begin{enumerate}
\item [{\rm (1)}]\ $T\leq N$.
  \item [{\rm (2)}]\ $\Omega_{0}, \Omega_{1}, \Omega'_{0}, \Omega'_{1}$ are four orbits of $T$.
  \item [{\rm (3)}]\ $\G[\Omega_{0}\cup \Omega_{1}\cup \Omega'_{0}\cup\Omega'_{1}]$ is a cycle of length $4t$.
  \item [{\rm (4)}]\ $\mb$ is a $G$-invariant partition of $V(\G)$.
\end{enumerate}
\end{lem}

\f\demo\ By Lemma~\ref{lem-not-block-4}, we see that $\R(a^l)\in N$, and so $T\leq N$. (1) holds. Since $T=\lg \R(a^l)\rg$ is assumed to be of order $t$, one has $T=\lg \R(a^{n/t})\rg$, and then one can obtain (2). By the adjacency rule of bi-Cayley graph, we can obtain (3).

Set $\Omega=\Omega_{0}\cup \Omega_{1}\cup \Omega'_{0}\cup\Omega'_{1}$ and $B=\Omega_{0}\cup\Omega_{1}$. By Lemma~\ref{lem-not-block-2}, $\G[\Delta]$ is a null graph, and so $B=\Delta\cap\Omega$. Since $\G$ has valency $3$, it follows that $\Delta\cup\Delta'$ is a block of imprimitivity of $G$ on $V(\G)$, and hence $\Omega$ is also a block of imprimitivity of $G$ on $V(\G)$ since $\G[\Omega]$ is a component of $\G[\Delta\cup\Delta']$. Since $\Delta$ is also a block of imprimitivity of $G$ on $V(\G)$, $B (=\Delta\cap\Omega)$ is a block of imprimitivity of $G$ on $V(\G)$.
Then $\mb=\{B^{\R(h)}\ |\ h\in H\}$ is a $G$-invariant partition of $V(\G)$.\hfill\qed

\begin{lem}\label{lem-not-block-6}
$T<N$ and the quotient graph $\G_{\mb}$ of $\G$ relative to $\mb$ is isomorphic to the cross ladder graph $\cl_{\frac{4n}{2t}}$. Moreover, $T$ is the kernel of $G$ acting on $\mb$.
\end{lem}

\f\demo
Let $K_{\mb}$ be the kernel of $G$ acting on $\mb$. Clearly, $T\leq K_\mb$. Let $B'=\Omega'_{0}\cup\Omega'_{1}$. Then $B'=B^{\R(b)}\in\mb$. Let $B^{\R(h)}\in\mb$ be adjacent to $B$ and $B^{\R(h)}\neq B'$.

Suppose that $\G[B\cup B^{\R(h)}]$ is a perfect matching. Since $G$ is transitive on $\mb$, $\G_{\mb}$ is a cycle of length $\frac{2n}{t}$. Clearly, $G/K_{\mb}$ is vertex-transitive but not edge-transitive on $\G_{\mb}$, so $G/K_{\mb}\cong D_{2n/t}$.
If $t=1$, then it is easy to see that $\G\cong \cl_{4n}$ which is a Cayley graph by Theorem~\ref{not-at-dihedrant}, a contradiction. If $t>1$, then since $\G[\Omega]=\G[B\cup B']$ is a cycle of length $4t$, $K_{\mb}$ acts faithfully on $B$, and so $K_{\mb}\leq\Aut(\G[B\cup B'])\cong D_{8t}$. Since $K_{\mb}$ fixes $B$, one has $|K_\mb|\mid 4t$, implying that $|G|=|K_\mb|\cdot\frac{2n}{t}\mid 8n$. As $|R(H)|=2n$ and $\R(H)$ is non-normal in $G$, one has $|K_\mb|=4t$ due to $T\leq K_\mb$. In view of the fact that $K_{\mb}\leq D_{8t}$, $K_{\mb}$ has a characteristic cyclic subgroup, say $J$, of order $2t$. Then we have $J\unlhd G$ because $K_{\mb}\unlhd G$. Clearly, $J$ is regular on $B$ and $J\cap N=T$, so $J\R(H)$ is regular on $V(\G)$. It follows from Proposition~\ref{IsCay} that $\G$ is a Cayley graph, a contradiction.

Therefore, $\G[B\cup B^{\R(h)}]$ is not a perfect matching. If $N=T$, then $B=\Delta$ and $B'=\Delta'$ are orbits of $M$, and then $\G[B\cup B^{\R(h)}]$ will be a perfect matching, a contradiction. Thus, $N>T.$

Now we are going to prove that $\G_{\mb}\cong \cl_{\frac{n}{2t}}$. Since $B$ is adjacent to $B^{\R(h)}$, $\Omega_{i}$ is adjacent to $\Omega_{j}^{\R(h)}$ for some $i,j\in \{0,1\}$. Then because $\Omega_{i}$ and $\Omega_{j}^{\R(h)}$ are orbits of $T$, $\G[\Omega_{i}\cup\Omega_{j}^{\R(h)}]$ is a perfect matching. This implies that $\G_{\mb}$ is of valency $3$, and so $K_{\mb}$ is intransitive on $B$. As every $B^{h}\in \mb$ is a union of two orbits of $T$ on $V(\G)$, $K_{\mb}$ fixes every orbit of $T$.
{ Since $N$ is cyclic, the normality of $N$ in $G$ implies that $T\unlhd G$. Clearly, $\Omega_{0}$ is adjacent to three pair-wise different orbits of $T$, so the quotient graph $\G_{T}$ of $\G$ relative to $T$ is of valency $3$.
Consequently, the kernel of $G$ acting on $V(\G_{T})$ is $T$.}
Then $K_{\mb}=T$.
Now $\R(H)/T\cong D_{2n/t}$ is regular on $\mb$, and so $\G_{\mb}$ is a Cayley graph over $\R(H)/T$. Furthermore, $G/T$ is not arc-transitive on $\G_{\mb}$. Since $\R(H)/T$ is non-normal in $G/T$, $\G_{\mb}$ is a non-normal Cayley graph over $\R(H)/T$. If $\G_{\mb}$ is arc-transitive, then by \cite[Theorem~1]{MP}, either $|\Aut(\G_{\mb})|=3k|\R(H)/T|$ with $k\leq 2$, or $\G_\mb$ has order $2\cdot p$ with $p=3$ or $7$. For the former, since $G/T$ is not arc-transitive on $\G_{\mb}$, one has $|G/T: \R(H)/T|\leq 2$, implying $\R(H)\unlhd G$, a contradiction. For the latter, we have $\frac{2n}{t}=6$ or $14$, implying $\frac{n}{t}=3$ or $7$.  It follows that $T$ is a maximal subgroup of $\lg \R(a)\rg$, and so $T=N$, a contradiction. Therefore, $\G_{\mb}$ is not arc-transitive. Since $\R(H)/T$ is non-normal in $G/T$, by Theorem~\ref{not-at-dihedrant}, one has $\G_{\mb}\cong \cl_{\frac{4n}{2t}}$, as required.\hfill\qed


\f{\bf Proof of Theorem~\ref{not-block}}\
By Lemma~\ref{lem-not-block-6}, we have $\G_{\mb}\cong \cl_{\frac{4n}{2t}}$. By the definition of $\cl_{\frac{4n}{2t}}$, we may partition the vertex set of $\G_{\mb}$ in the following way:
\[V(\G_{\mb})=V_0\cup V_1\cup \cdots V_{\frac{2n}{2t}-2}\cup V_{\frac{2n}{2t}-1},\ {\rm where}\ V_i=\{B_i^0, B_i^1\}, i\in\mz_{\frac{2n}{2t}}\]
and
\[E(\G_\mb)=\{\{B_{2i}^r, B_{2i+1}^r\}, \{B^r_{2i+1}, B^s_{2i+2}\}\mid i\in\mz_{\frac{n}{2t}}, r,s\in\mz_2\}.\]

Assume that $B_0^0=B$ and $B_1^0=B'$. Recall that $B=\Omega_{0}\cup\Omega_{1}$ and $B'=\Omega'_{0}\cup\Omega'_{1}=B^{\R(b)}$.
Moreover, $\Omega_{0}, \Omega_{1}, \Omega'_{0}$ and $\Omega'_{1}$ are four orbits of $T$.
Then every $B_i^j\in\mb$ is just a union of two orbits of $T$. For convenience, we may let
\[B_i^j=\Omega_{i0}^j\cup\Omega_{i1}^j, i\in\mz_{\frac{2n}{2t}}, j\in\mz_2,\]
where $\Omega_{i0}^j, \Omega_{i1}^j$ are two orbits of $T$.
For $B=B_0^0$, we let $\Omega_{0}=\Omega_{00}^0$ and $\Omega_{1}=\Omega_{01}^0$, and for $B'=B_1^0$, we let
$\Omega_{0}'=\Omega_{10}^0$ and $\Omega_{1}'=\Omega_{11}^0$.

For convenience, in the remainder of the proof, we shall use $\C_{4t}$ to denote a cycle of length $4t$, and we also call $\C_{4t}$ a $4t$-cycle.
Recall that $\G[B\cup B']=\G[B_0^0\cup B_1^0]\cong \C_{4t}$, and that the edges between $\Omega_{0i}^0(=\Omega_{i})$ and $\Omega_{1j}^0(=\Omega_{j}')$ form a perfect matching for all $i,j\in\mz_2$.
Since $T\unlhd G$, the quotient graph $\G_T$ of $\G$ relative to $T$ has valency $3$.
So the edges between any two adjacent orbits of $T$ form a perfect matching.

\begin{figure}[htbp]
\begin{center}
\unitlength 4mm
\begin{picture}(25,9)

\put(11, 8){\line(1, 0){4}}\put(11, 8){\line(2, -1){4}}
\put(11, 6){\line(1, 0){4}}\put(11, 6){\line(2, 1){4}}

\put(11, 3){\line(1, 0){4}}\put(11, 3){\line(2, -1){4}}
\put(11, 1){\line(1, 0){4}}\put(11, 1){\line(2, 1){4}} 

\put(11, 8){\line(-1, 0){3}}\put(15, 8){\line(1, 0){3}}
\put(11, 6){\line(-1, -1){3}}\put(15, 6){\line(1, -1){3}}

\put(11, 3){\line(-1, 1){3}}\put(15, 3){\line(1, 1){3}}
\put(11, 1){\line(-1, 0){3}}\put(15, 1){\line(1, 0){3}}

\put(10.7, 7.7){$\circ$}\put(10.7, 5.7){$\circ$}
\put(10.7, 2.7){$\circ$}\put(10.7, 0.7){$\circ$}

\put(14.7, 7.7){$\circ$}\put(14.7, 5.7){$\circ$}
\put(14.7, 2.7){$\circ$}\put(14.7, 0.7){$\circ$}

\put(7.7, 7.7){$\circ$}\put(7.7, 5.7){$\circ$}
\put(7.7, 2.7){$\circ$}\put(7.7, 0.7){$\circ$}

\put(17.7, 7.7){$\circ$}\put(17.7, 5.7){$\circ$}
\put(17.7, 2.7){$\circ$}\put(17.7, 0.7){$\circ$}

\put(11,6.9){\oval(1,3.2)}\put(11,2.2){\oval(1,3.2)}
\put(8,6.9){\oval(1,3.2)}\put(8,2.2){\oval(1,3.2)}

\put(15,6.9){\oval(1,3.2)}\put(15,2.2){\oval(1,3.2)}
\put(18,6.9){\oval(1,3.2)}\put(18,2.2){\oval(1,3.2)}

\put(5.5, 6.7){$\cdots$}\put(5.5, 1.7){$\cdots$}
\put(19, 6.7){$\cdots$}\put(19, 1.7){$\cdots$}

\put(7.7, 8.7){{\footnotesize$B_{-1}^0$}}\put(10.8, 8.7){{\footnotesize$B_0^0$}}
\put(14.5, 8.7){{\footnotesize$B_1^0$}}\put(17.5, 8.7){{\footnotesize$B_2^0$}}
\put(6.7, -0.6){{\footnotesize$B_{-1}^1$}}\put(10.5, -0.6){{\footnotesize$B_0^1$}}
\put(14.5, -0.6){{\footnotesize$B_1^1$}}\put(17.6, -0.6){{\footnotesize$B_2^1$}}

\put(9.5, 7.3){{\footnotesize$\Omega_{00}^0$}}\put(9.5, 5.6){{\footnotesize$\Omega_{01}^0$}}
\put(15.5, 7.3){{\footnotesize$\Omega_{10}^{0}$}}\put(15.5, 5.6){{\footnotesize$\Omega_{11}^{0}$}}

\put(9.5, 2.8){{\footnotesize$\Omega_{00}^1$}}\put(9.5, 1.1){{\footnotesize$\Omega_{01}^1$}}
\put(15.5, 2.8){{\footnotesize$\Omega_{10}^{1}$}}\put(15.5, 1.1){{\footnotesize$\Omega_{11}^{1}$}}

\put(5, 2.8){{\footnotesize$\Omega_{(-1)0}^1$}}\put(5, 1.1){{\footnotesize$\Omega_{(-1)1}^1$}}
\put(5, 7.8){{\footnotesize$\Omega_{(-1)0}^0$}}\put(5, 5.6){{\footnotesize$\Omega_{(-1)1}^0$}}

\put(18.6, 2.8){{\footnotesize$\Omega_{20}^1$}}\put(18.6, 1.1){{\footnotesize$\Omega_{21}^1$}}
\put(18.6, 7.3){{\footnotesize$\Omega_{20}^0$}}\put(18.6, 5.6){{\footnotesize$\Omega_{21}^0$}}

\end{picture}
\end{center}\vspace{-.1cm}
\caption{The sketch graph of $\G_{\mb}$} \label{Fig-2.4}
\end{figure}
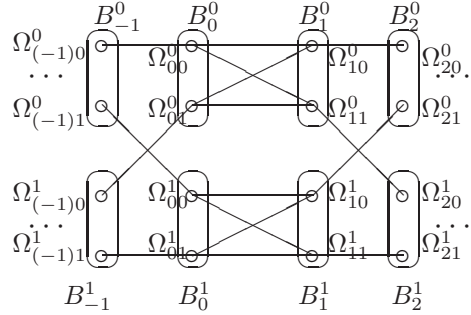

From the construction of $\G_\mb,$ one may see that there exists $g\in G$ such that
$\{V_{0}, V_1\}^g=\{V_{2i}, V_{2i+1}\}$ for each $i\in\mz_{\frac{n}{2t}}$.
 So for each $i\in\mz_{\frac{n}{2t}}, r\in\mz_2$, we may assume that
$\G[B_{2i}^r\cup B_{2i+1}^r]\cong \C_{4t}$, and $\Omega_{(2i)s}^r\sim\Omega_{(2i+1)t}^r$ for all $s,t\in\mz_2$. (Here $\Omega_{(2i)s}^r\sim\Omega_{(2i+1)t}^r$ means that $\Omega_{(2i)s}^r$ and $\Omega_{(2i+1)t}^r$ are adjacent in $\G_\mb$.)
Again, from the construction of $\G_\mb$, we may assume that
\[\Omega_{(2i+2)0}^0\sim\Omega_{(2i+1)0}^0, \Omega_{(2i+2)1}^0\sim\Omega_{(2i+1)0}^1, \Omega_{(2i+2)1}^1\sim\Omega_{(2i+1)0}^1, \Omega_{(2i+2)1}^1\sim\Omega_{(2i+1)0}^0,\] for each $i\in\mz_{\frac{n}{2t}}$. We draw a local subgraph of $\G_{\mb}$ in Figure~\ref{Fig-2.4}.
Observing that every $V_i=\{B_i^0, B_i^1\}$ with $i\in\mz_{\frac{2n}{2t}}$ is a block of imprimitivity of $G/K_\mb$ acting on $V(\G_\mb)$.
So every $B_i^0\cup B_i^1$  with $i\in\mz_{\frac{2n}{2t}}$ is a block of imprimitivity of $G$ acting on $V(\G)$.
Let $E$ be the kernel of $G$ acting on the block system $\Lambda=\{B_i^0\cup B_i^1\mid i\in\mz_{\frac{2n}{2t}}\}$.
Then $G/E\cong D_{\frac{n}{t}}$ acts regularly on $\Lambda$. Clearly, $\R(H)$ is also transitive on $\Omega$, so $G/E=\R(H)E/E$.
By Lemma~\ref{lem-not-block-6}, $T$ is a the kernel of $G$ acting on $\mb$. So $E/T$ is an elementary $2$-group.
%
From $\R(H)/(\R(H)\cap E)\cong D_{\frac{n}{t}}$ it follows that $\R(H)\cap E=\lg\R(a^{\frac{n}{2t}})\rg\cong\mz_{2t}$, and so $(\R(H)\cap E)/T$ is a normal subgroup of $G/T$ of order $2$. This implies that $B_i^1=(B_i^0)^{\R(a^{\frac{n}{2t}})}$ for $i\in\mz_{\frac{2n}{2t}}$. We may further assume that $\Omega_{01}^1=(\Omega_{00}^0)^{\R(a^{\frac{n}{2t}})}\subseteq B_0^1$.
So $\Omega_{00}^0\cup\Omega_{01}^1$ is just the orbit of $\lg\R(a^{\frac{n}{2t}})\rg$ containing $1_0$.

Observing that $\Omega_{10}^0\sim\Omega_{20}^0$ and the edges between them are of the form $\{g_0, (ba^ig)_0\}$ with $g_0\in \Omega_{10}^0$, one has $\Omega_{20}^0=ba^i\Omega_{10}^0=ba^i(\Omega_{00}^0)^{\R(b)}=(\Omega_{00}^0)^{\R(a^{-i})}$.
So $\Omega_{20}^1\subseteq (B_0^1)^{\R(a^{-i})}$.

Since $B_1^0=B'=B^{\R(b)}=(B_{0}^0)^{\R(b)}$, one has $B_1^1=(B_{0}^1)^{\R(b)}$. Recall that $1_1\in \Omega_{11}^0=\Omega_1'$ and $1_1$ is adjacent to $1_0\in\Omega_{00}^0=\Omega_0$ and $(ba^l)_1\in \Omega_{01}^0=\Omega_1$. As we assume that $\Omega_{11}^0\sim \Omega_{20}^1$, $1_1$ is adjacent to some vertex in $\Omega_{20}^1$. So $\Omega_{20}^1\subseteq H_1$ and hence \[\Omega_{20}^1=(\Omega_{00}^1)^{\R(a^{-i})}=(\Omega_{01}^0)^{\R(a^{\frac{n}{2t}})\R(a^{-i})}=(\Omega_{01}^0)^{\R(a^{\frac{n}{2t}-i})}=\{(ba^{k\frac{n}{t}})_{1}\ |\ 0\leq k\leq t-1\}^{\R(a^{\frac{n}{2t}-i})}.\]
So we have the following claim.

\medskip
\f{\bf Claim~1}\ $L=\{ba^l, ba^{k\frac{n}{t}+\frac{n}{2t}-i}\}$ and $R=\{ba^i, b\}$, where $|\R(a^l)|=t$, $i\in\mz_n$ and $0\leq k\leq t-1$. \medskip

Let $G_{1_0}^*$ be the kernel of $G_{1_0}$ acting on the neighborhood of $1_0$ in $\G$.
Then $G_{1_0}^*\leq E_{1_0}$. Recall that for each $i\in\mz_{\frac{n}{2t}}, r\in\mz_2$,
$\G[B_{2i}^r\cup B_{2i+1}^r]\cong \C_{4t}$ and the edges between $B_{2i+1}^0\cup B_{2i+1}^1$ and $B_{2i+2}^0\cup B_{2i+2}^1$ form a perfect matching.
It follows that $E$ acts faithfully on each $B_i^0\cup B_i^1$. Clearly, $G_{1_0}^*\leq E_{1_0}$, so $G_{1_0}^*$ acts faithfully on each $B_i^0\cup B_i^1$.

\medskip

\f{\bf Claim~2}\ If $t>2$ then $G_{1_0}^*=1$, and if $t=2$ then $G_{1_0}^*\leq\mz_2$ and $3\mid n$.

Assume that $t\geq 2$. Since $\G[B_{0}^0\cup B_{1}^0]\cong \C_{4t}$, $G_{1_0}^*$ fixes every vertex in $B_{0}^0$, and so fixes every vertex in $\Omega_{(-1)0}^0$ since $\Omega_{(-1)0}^0\sim\Omega_{00}^0$ (see Figure~\ref{Fig-2.4}). This implies that $G_{1_0}^*$ fixes $\Omega_{(-1)1}^0$ setwise, and so fixes $\Omega_{00}^1$ setwise since $\Omega_{(-1)1}^0\sim\Omega_{00}^1$. Consequently, $G_{1_0}^*$ also fixes $\Omega_{01}^1$ setwise. Similarly, by considering the edges between $B_1^0\cup B_1^1$ and $B_2^0\cup B_2^1$, we see that $G_{1_0}^*$ fixes both $\Omega_{10}^1$ and $\Omega_{11}^1$ setwise. Recall that the edges between $\Omega_{0i}^1$ and $\Omega_{1j}^1$ form a perfect matching for $i,j\in\mz_2$. As  $\G[B_{1}^0\cup B_{1}^1]\cong \C_{4t}$, $G_{1_0}^*$ acts faithfully on $\Omega_{00}^1$ (or $\Omega_{01}^1$), and so $G_{1_0}^*\leq\mz_2$.

If $t>2$, then since $\G[B_{-2}^0\cup B_{-1}^0]\cong \C_{4t}$, $G_{1_0}^*$ will fix every vertex in this cycle, and in particular,
$G_{1_0}^*$ will fix every vertex in $\Omega_{(-1)1}^0$. As $\Omega_{(-1)1}^0\sim\Omega_{00}^1$, $G_{1_0}^*$ will fix every vertex in $\Omega_{00}^1$. Since $G_{1_0}^*$ acts faithfully on $\Omega_{00}^1$, one has $G_{1_0}^*=1$.


Let $t=2$. We shall show that $3\mid n$. Then $T=\lg\R(a^{\frac{n}{2}})\rg$. Recall that $(\R(H)\cap E)/T $ is a  normal subgroup of $G/T$ of order $2$. Let $M=\R(H)\cap E$. Then $M$ is a normal subgroup of $G$ of order $4$. Since $\R(H)$ is dihedral, one has $M=\lg\R(a^\frac{n}{4})\rg$.
Let $C=C_G(M)$. Then $\R(a)\in C$ and $\R(b)\notin C$. It follows that $C$ is a proper subgroup of $G$. Since $G/E$ acts regularly on $\Lambda$, $C_{1_0}$ fixes every element in $\Lambda$. Since $C_{1_0}$ centralizes $M$, $C_{1_0}$ fixes every vertex in the orbit $\Omega_{00}^0\cup\Omega_{01}^1$ of $M$ containing $1_0$. Clearly, $C_{1_0}\leq G_{1_0}$, so $C_{1_0}/(C_{1_0}\cap G_{1_0}^*)\leq\mz_2$. As we have shown that $G_{1_0}^*$ acts faithfully on $\Omega_{01}^1$, it follows that $C_{1_0}\cap G_{1_0}^*=1$ since $C_{1_0}$ fixes $\Omega_{01}^1$ pointwise, and hence $C_{1_0}\leq\mz_2$. On the other hand, as $G_{1_0}^*\leq\mz_2$, one has $|G|\mid 4\cdot 4n=16n$.
Since $C<G$ and $\R(a)\in C$, one has $|C|=kn$ with $k\mid 8$.

Suppose that $3\nmid n$. For any odd prime divisor $p$ of $n$, let $P$ be a Sylow $p$-subgroup of $\lg\R(a)\rg$. Then $P$ is also a Sylow $p$-subgroup of $C$. If $P$ is not normal in $C$, then by Sylow's theorem, we have $|C: N_C(P)|=k'p+1\mid 8$ for some integer $k'$. Since $p\neq 3$, one has $p=7$ and $k'=1$. This implies that $|C|=8|N_C(P)|$, and so $|C|=8n$ due to $\R(a)\in C$ and $C<G$. Since $C_{1_0}\leq\mz_2$, one has $|C: C_{1_0}|\geq 4n$, and so $C$ is transitive on $V(\G)$. Moreover, we have $C_C(P)=N_C(P)=\lg\R(a)\rg$. By Burnside theorem, $C$ has a normal subgroup $M$ such that $C=M\rtimes P$. Then the quotient graph $\G_M$ of $\G$ relative to $M$ would be a cycle of length $|P|$, and the subgraph induced by each orbit of $M$ is just a perfect matching. This implies that $M$ is just the kernel of $G$ acting on $V(\G_M)$. Furthermore, $C/M$ is a vertex-transitive subgroup of $\Aut(\G_M)$. Since $\G_M$ is a cycle, $C/M$ must contain a subgroup, say $B/M$ acting regularly on $V(\G_M)$. Then $B$ will be regular on $V(\G)$, and so by Proposition~\ref{IsCay}, $\G$ is a Cayley graph, a contradiction. Therefore, $P\unlhd C$, and since $C\unlhd G$, one has $P\unlhd G$, implying $P\leq N$. By the arbitrariness of $P$, $n/|N|$ must be even, contrary to Lemma~\ref{lem-not-block-1}.
Thus, $3\mid n$, as claimed. \medskip

The following claim shows that $t=1$ or $2$.

\medskip
\f{\bf Claim~3}\ $t\leq 2$.

By way of contradiction, suppose that $t>2$. Let $C=C_G(T)$. Then $\lg\R(a)\rg\leq C$ and $\R(H)\nleq C$ since $|T|=t>2$. Clearly, $C_{1_0}\leq E_{1_0}$. As $C_{1_0}$ centralizes $T$, $C_{1_0}$ will fixes every vertex in $\Omega_{00}^0$ since $\Omega_{00}^0$ is an orbit of $T$ containing $1_0$. Since $\G[B_{0}^0\cup B_{1}^0]\cong \C_{4t}$, $C_{1_0}$ fixes every vertex in this $4t$-cycle, and so $C_{1_0}\leq G_{1_0}^*=1$ (by Claim~2). Thus, $C$ acts semiregularly on $V(\G)$.
If $C=\lg \R(a)\rg$, then by N/C-theorem, we have $G/\lg \R(a)\rg=G/C\leq\Aut(T)$. Since $T\leq N\leq \lg \R(a)\rg$ is cyclic, $\Aut(T)$ is abelian. It then follows that $\R(H)/C\unlhd G/C$, and hence $\R(H)\unlhd G$, a contradiction.
If $C>\lg \R(a)\rg$, then $|C|=2n$ because $\G$ is non-Cayley. Since $H_0$ and $H_1$ are not blocks of imprimitivity of $G$ on $V(\G)$, $C$ does not fix $H_0$ setwise, and so $\R(H)C$ is transitive on $V(\G)$. Clearly, $\R(H)\cap C=\lg\R(a)\rg$, so $|\R(H)C|=|\R(H)||C|/|\lg\R(a)\rg|=4n$. It follows that $\R(H)C$ is regular on $V(\G)$, contradicting that $\G$ is non-Cayley.\medskip

By Claim~3, we only need to consider the following two cases:\medskip

\f{\bf Case~1}\ $t=1$.

In this case, by Claim~1, we have $R=\{b, ba^i\}$ and $L=\{b, ba^{\frac{n}{2}-i}\}$. For convenience, we let $n=2\ell$. Then
$R=\{b, ba^i\}$ and $L=\{b, ba^{\ell-i}\}$.

By Proposition~\ref{properties}~(1), the connectedness of $\G$ implies that $\lg a^{i}, a^{\ell}\rg=\lg a\rg$. Then either $(i, 2\ell)=1$, or $i=2k$ with $(k, 2\ell)=1$ and $\ell$ is odd.
Recall that $H=\lg a, b\mid a^{2\ell}=b^2=1, bab=a^{-1}\rg$. For any $\ld\in\mz_{2\ell}^*$, let $\a_\ld$ be the automorphism of $H$ induced by the map
\[a^{\ld}\mapsto a,~b\mapsto b.\]
So if $(i, 2\ell)=1$, then we have \[(R, L)^{\a_i}=(\{b, ba\}, \{b, ba^{\ell-1}\}),\]
and if $i=2k$ with $(k, 2\ell)=2$ and $\ell$ is odd,
then we have \[(R, L)^{\a_k}=(\{b, ba^2\}, \{b, ba^{\ell-2}\}).\]
So by Proposition~\ref{properties}~(3), we have
\[(R, L, S)\equiv(\{b, ba\}, \{b, ba^{\ell-1}\}, \{1\})\ {\rm or}\ (\{b, ba^2\}, \{b, ba^{\ell-2}\}, \{1\}) (\ell\ {\rm is\ odd}).\]

Suppose that $\ell$ is even. Then $(R, L, S)\equiv(\{b, ba\}, \{b, ba^{\ell-1}\}, \{1\})$.
Since $\ell$ is even, one has $(2\ell, \ell+1)=1$ and $(\ell+1)^2\equiv 1~(\mod 2\ell)$.
Then it is easy to check that $\a_{\ell+1}$ is an automorphism of $H$ of order $2$ that swaps $\{b, ba\}$ and $\{b, ba^{\ell-1}\}$.
By Proposition~\ref{normalizer}, we have $\d_{\a_{\ell+1},1,1}\in I$, and then $\G\cong\BiCay(H, \{b,ba\},\{b,ba^{\ell+1}\},\{1\})$ is a Cayley graph, a contradiction.

Now we assume that $n=2\ell$ with $\ell=2m+1$ for some integer $m$.
Let
\[\G_{1}=\BiCay(H, \{b, ba\}, \{b, ba^{2m}\}, \{1\}), \G_{2}=\BiCay(H, \{b, ba^{2}\}, \{b, ba^{2m-1}\}, \{1\}).\]
Direct calculation shows that $(n, 2m-1)=1$, and $2m(2m-1)\equiv2\ (\mod n)$.
Then the automorphism $\a_{2m-1}: a\mapsto a^{2m-1},~b\mapsto b$ maps
the pair of two subsets $(\{b, ba\}, \{b, ba^{2m}\})$ to $(\{b, ba^{2m-1}\}, \{b, ba^{2}\}).$
So, we have $(R, L, S)\equiv(\{b, ba\}, \{b, ba^{2m}\}, \{1\})$.
By Lemma~\ref{EX(m,2)} and Theorem~\ref{EX(m,2)-Cayley}, $\G\cong \mcl(4m, 2)$ and $\G$ is non-Cayley if and only if $3\nmid (2m+1)$.
Note that $3\nmid (2m+1)$ is equivalent to $m\not\equiv1~(\mod 3)$. So we obtain the first family of graphs in Theorem~\ref{not-block}.\medskip

\f{\bf Case~2}\ $t=2$.

In this case, by Claim~1, we have $R=\{b, ba^i\}$ and $L=\{ba^{\frac{n}{2}}, ba^{\frac{3n}{4}-i}\}$ or $\{ba^{\frac{n}{2}}, ba^{\frac{n}{4}-i}\}$.
We still use the following notation: For any $\ld\in\mz_{2\ell}^*$, let $\a_\ld$ be the automorphism of $H$ induced by the map
\[a^{\ld}\mapsto a,~b\mapsto b.\]
Note that \[(\{b, ba^i\}, \{ba^{\frac{n}{2}}, ba^{\frac{3n}{4}-i}\})^{\a_{-1}}=(\{b, ba^{-i}\}, \{ba^{\frac{n}{2}}, ba^{\frac{n}{4}-(-i)}\}).\]
By replacing $-i$ by $i$, we may always assume that
\[(R, L)=(\{b, ba^i\}, \{ba^{\frac{n}{2}}, ba^{\frac{n}{4}-i}\}).\]

By Claim~2, we have $3\mid n$. So we may assume that $n=12m$ for some integer $m$. Then we have
\[(R, L)=(\{b, ba^i\}, \{ba^{6m}, ba^{3m-i}\}).\]
Since $\G$ is connected, by Proposition~\ref{properties}, we have $\lg a^i, a^{3m}\rg=\lg a\rg$.
If $m$ is odd, by Lemma~\ref{Cayley}, $\G$ will be a Cayley graph which is impossible. Thus, $m$ is even.
It then follows that $\lg a^i\rg\cap \lg a^{3m}\rg>1$ since $\lg a^i, a^{3m}\rg=\lg a^i\rg\lg a^{3m}\rg=\lg a\rg$.
Since $\lg a^{3m}\rg\cong\mz_4$, one has $|\lg a^i\rg\cap \lg a^{3m}\rg|=2$ or $4$. For the former, we would have
$|\lg a^i\rg|=6m$, and since $m$ is even, one has $4\mid |\lg a^i\rg|$, and hence $a^{3m}\in \lg a^i\rg$, a contradiction.
Thus, we have $|\lg a^i\rg\cap \lg a^{3m}\rg|=4$, that is, $\lg a^i\rg=\lg a\rg.$
So $(i, 12m)=1$, and then $\a_i\in\Aut(H)$ which maps $(\{b, ba^{i}\}, \{ba^{6m},ba^{3m-i}\})$ to
$(\{b, ba\}, \{ba^{6m}, ba^{3m-1}\})$ or $(\{b, ba\}, \{ba^{6m}, ba^{-3m-1}\})$.
Then \[(R, L, S)\equiv (\{b, ba\}, \{ba^{6m}, ba^{3m-1}\}, \{1\})\ {\rm or}\ (\{b, ba\}, \{ba^{6m}, ba^{-3m-1}\}), \{1\}.\]
If $m\equiv 2\ (\mod 4)$, then by Lemma~\ref{Cayley2}, we see that $\G$ will be a Cayley graph, a contradiction.
Thus, $m\equiv 0\ (\mod 4)$. Clearly, $(3m-1, 12m)=1$, and hence the map $a\mapsto a^{3m-1}, b\mapsto ba^{6m}$ induces an automorphism, say $\b$ of $H$. It is easy to check that
$$
(\{b,ba\}, \{ba^{6m}, ba^{3m-1}\})^{\b}=(\{ba^{6m},ba^{-3m-1}\},\{b,ba\}).$$
Thus, \[(R, L, S)\equiv (\{b, ba\}, \{ba^{6m}, ba^{3m-1}\}, \{1\}).\]
By Proposition~\ref{non-Cayley-2}, $\G$ is a non-Cayley graph. Let $m=4\ell$ for some integer $\ell$. Then $n=48\ell$ and then we get the second family of graphs in Theorem~\ref{not-block}. This completes the proof of Theorem~\ref{not-block}.
\hfill\qed

\subsection{Proof of Theorem~\ref{all-non-cay}}

By \cite[Theorem 1.2]{ZZ}, if $\G$ is $0$- or $1$-type, then $\G$ is a Cayley graph. Let $\G$ be of $2$-type. Suppose that $\G$ is a non-Cayley graph. Let $G\leq\Aut(\G)$ be minimal subject to that $\R(H)\leq G$ and $G$ is transitive on $V(\G)$. If $H_0$ and $H_1$ are blocks of imprimitivity of $G$ on $V(\G)$, then by Proposition~\ref{block}, we obtain the first two families of graphs of Theorem~\ref{all-non-cay}. Otherwise, $H_0$ and $H_1$ are not blocks of imprimitivity of $G$ on $V(\G)$, by Theorem~\ref{not-block}, we obtain the last two families of graphs of Theorem~\ref{all-non-cay}. \hfill\qed

\medskip
\f {\bf Acknowledgements:} This work was supported by the National
Natural Science Foundation of China (11671030).

\end{document}